\documentclass[Journal]{IEEEtran}
\usepackage{amssymb}
\usepackage{graphicx,subfigure}
\usepackage{latexsym,amsmath,amsthm,amssymb,amsfonts}
\usepackage{mathrsfs}
\usepackage{amscd}
\usepackage[ansinew]{inputenc}
\usepackage{color}
\usepackage{cite}
\usepackage[all]{xy}
\usepackage{balance}
\usepackage[justification=centering]{caption}
\usepackage{makecell}
\newsavebox{\mybox}
\newcolumntype{X}[1]{%
>{\begin{lrbox}{\mybox}}%
c%
<{\end{lrbox}\makecell[#1]{\usebox\mybox}}%
}
\usepackage{makecell}

\def\QED{\mbox{\rule[0pt]{1.5ex}{1.5ex}}}
\def\endproof{\hspace*{\fill}~\QED\par\endtrivlist\unskip}

\newtheorem{theorem}{Theorem}

\newtheorem{remark}{Remark}
\newtheorem{lemma}{Lemma}
\newtheorem{assumption}{Assumption}

\renewcommand{\baselinestretch}{1}

\begin{document}

\title {Distributed Continuous-Time and Discrete-Time Optimization With Nonuniform Unbounded Convex Constraint Sets and Nonuniform Stepsizes
}
\author{\mbox{Peng Lin}\thanks{Peng Lin, Chunhua Yang and Weihua Gui are with the School of Information Science and Engineering, Central South University, Changsha, China. Wei Ren is with the Department of Electrical and Computer Engineering, University of California, Riverside, USA.
E-mail: lin$\_$peng0103@sohu.com,
ren@ee.ucr.edu, ychh@csu.edu.cn, gwh@csu.edu.cn.
This work was supported by the National Science Foundation under Grant ECCS-1307678 and ECCS-1611423, the Foundation for Innovative Research Groups of
the National Natural Science Foundation of China (61321003), the 111 Project
(B17048), the National Natural Science Foundation of China (61573082, 61203080),
and the Innovation-driven Plan
at Central South University.}, Wei Ren, Chunhua Yang and Weihua Gui}

    \markboth{}
        {}
\maketitle

\begin{abstract}                         
This paper is devoted to distributed continuous-time and discrete-time optimization problems with nonuniform convex constraint sets and nonuniform stepsizes
for general differentiable convex objective functions. The communication graphs are not required to be strongly connected at any time, the gradients of the local objective functions are not required to be bounded when their independent variables tend to infinity, and the constraint sets are not required to be bounded. For continuous-time multi-agent systems, a
 distributed continuous
algorithm is first introduced where the stepsizes and the convex constraint sets are both nonuniform. It is shown that all agents reach a consensus while minimizing the team objective function even when the constraint sets are unbounded.
After that, the obtained results are extended to discrete-time multi-agent systems and then the case where each agent remains in a corresponding convex constraint set is studied. To ensure all agents to remain in a bounded region, a switching mechanism is introduced in the algorithms. It is shown that the distributed optimization problems can be solved, even though the discretization of the algorithms might deviate the convergence of the agents from the minimum of the objective functions. Finally, numerical examples are included to show the obtained
theoretical results.
\vspace{0.1cm}

\noindent{\bf Keywords}: Distributed Optimization, Nonuniform Step-Sizes, Nonuniform Convex Constraint Sets
\end{abstract}

\section{Introduction}
As an important research direction of control theory,
distributed optimization problems {{for}} multi-agent systems have attracted
more and more attention from the control community
\cite{angelia,angelia1,nedd,Elia,cotes,Cortes3,Wei,srivast,nedic0,Zhu,Johansson,lu,lup,shi,liu,Kvaternik,Liuliu,Zavlanos,N1,yuan,Lin2014,hong3,hong4,lei,ren, hu,hong2,Ch,linren3,xie6, linren4}.
The goal of a distributed optimization problem for a multi-agent system is to minimize a desired team objective function cooperatively in a distributed way
where each agent can only have access to partial information of the team objective function. During the past few years,
several results have been obtained for distributed optimization problems.
For example, {{article}} \cite{angelia} introduced a discrete-time projection algorithm for multi-agent systems with state constraints
 and proved that the optimization problems can be solved when the communication topologies are jointly strongly connected and balanced. Articles \cite{shi} and \cite{liu}  studied  a continuous-time
version of the work in \cite{angelia} with convex constraint sets. Article \cite{Elia} gave a distributed continuous-time dynamic  algorithm for distributed optimization, and subsequently, on this basis, {{articles}} \cite{cotes,Cortes3} studied the distributed optimization problem for general strongly connected balanced directed graphs and gave the estimate of the convergence rate of the algorithm. Other works about distributed optimization problems could be found in articles \cite{angelia1,nedd,Wei,srivast,nedic0,Zhu,Johansson,lu,lup,Kvaternik,Liuliu,Zavlanos,N1,yuan,Lin2014,hong3,hong4,lei,ren,hu,hong2,Ch,linren3,xie6,linren4} and the references therein, where new algorithms, e.g., distributed Newton, approximate dual subgradient and zero-gradient-sum algorithms, were given or more complicated cases, e.g., second-order dynamics, time-varying or nonconvex functions, fixed or asynchronous stepsizes and noise, were considered.

Though many excellent results have been obtained for the distributed optimization problem, many issues need be further studied, e.g., general convex functions, nonuniform convex constraints and nonuniform stepsizes. For the issues of general convex functions and nonuniform convex constraints, most of the existing results require the gradients or subgradients of the convex functions to be bounded and the convex constraints to be identical and little attention has been paid to general convex functions and nonuniform convex constraints, in particular for multi-agent systems with general directed balanced graphs and unbounded gradients.
 For example, articles \cite{cotes,Cortes3}
 studied general convex functions but assumed them to be globally Lipschitz and the graphs are assumed to be strongly connected and balanced. Article \cite{liu} studied coercive convex functions with unbounded subgradients but the results are limited to the continuous-time multi-agent systems and the convex constraints sets are assumed to be identical for all agents.
Article \cite{angelia} studied nonuniform convex constraints but the communication graph is complete and all the edge weights are assumed to be equal.
{Founded on \cite{angelia}, articles \cite{hong4,lei} gave some results on nonuniform convex constraints but the communication graph is constant and connected, and the objective functions are assumed to be strongly convex or some intermediate variables need be transmitted besides the agent states.}
 Article
\cite{linren3} studied a distributed optimization problem with nonuniform convex constraints and gave conditions to guarantee the optimal convergence of the team objective function, but the subgradients and the convex constraint sets are both bounded. Note that \cite{angelia,liu, hong4,lei,linren3} all adopt a uniform stepsize. For the issue of nonuniform stepsizes, currently, there are few works concerned about this issue. {Articles \cite{srivast,nedic0} studied the distributed optimization problem with nonuniform stepsizes in a stochastic setting, where the communication graphs are required to be undirected and connected.} Article \cite{linren4} introduced a kind of nonuniform stepsizes but the discontinuous algorithms were employed to realize the consensus of all agents and the convex constraint sets are assumed to be identical. {Article \cite{xie6} also studied the distributed optimization problem with nonuniform stepsizes but some intermediate variables need be transmitted besides the agent states in order to track the average of the gradients.}
{\footnotesize{\begin{table}[!hbp]
\begin{tabular}{|X{cc}|X{cc}|X{cc}|}
\hline
\scriptsize{distributed optimization} &\makecell[tl]{\scriptsize{nonuniform stepsizes}\\\scriptsize{~~~~~~~~~~(A)} }& \makecell[tl]{\scriptsize{nonuniform stepsizes}\\{\scriptsize~~~~~~~~~~(B)}}\\
\hline
 \makecell[tl]{\scriptsize{nonuniform convex}\\\scriptsize{constraints}} &  \centering\scriptsize{unsolved} & \scriptsize{unsolved} \\
\hline
\makecell[tl]{\scriptsize{jointly strongly}\\\scriptsize{connected graphs}}&  \centering\scriptsize{unsolved} & \scriptsize{unsolved} \\
\hline
 \centering{\scriptsize{unbounded gradients}}&\makecell[tl]{\scriptsize{partly solved with certain}\\\scriptsize{ assumptions, e.g. \cite{linren4,xie6}}}&\makecell[tl]{\scriptsize{partly solved in \cite{linren4} with}\\\scriptsize{discontinuous algorithms}}\\
\hline
\end{tabular}
\caption{\upshape Existing results on distributed optimization with nonuniform stepsizes. }
\scriptsize{Nonuniform stepsizes (A) denotes the case where some intermediate variables need to be transmitted besides the agents' states while nonuniform stepsizes (B) denotes the case where no intermediate variables need to be transmitted. The main contributions of this paper lie in dealing with the above three different situations simultaneously using continuous algorithms in the case of nonuniform stepsizes (B).}
\end{table}}}

 In this paper, we are interested in studying distributed continuous-time and discrete-time optimization problems with nonuniform convex constraint sets and nonuniform stepsizes
for general differentiable convex objective functions. The communication graphs might not be strongly connected at any time and it is only required that  the union of the communication graphs among the time intervals of a certain length be strongly connected. The gradients of the local objective functions considered might not be bounded when their independent variables tend to infinity.
First, a distributed continuous-time
algorithm is introduced where the stepsizes and the convex constraint sets are both nonuniform. Nonuniform stepsizes mean that the weights of the gradients of the local objective functions in the control input of the agents are nonuniform. That is, the optimal convergence rates of the local objective functions are different, which has great possibility to result in the destruction of the optimal convergence of the team objective function.
 The existing works (e.g., \cite{angelia}) usually assumed the stepsizes are uniform and took no consideration of nonuniform stepsizes, and hence their approaches are hard to be applied for the case of nonuniform stepsizes.
 Our approach is to introduce a kind of stepsizes such that the stepsizes of each agent are constructed only based on its own states and the differences between the stepsizes of all agents vanish to zero as time evolves. Though the stepsizes of all agents tend to the same as time evolves, their differences still heavily affect the consensus stability and optimal convergence of the system, especially when the communication graphs are not strongly connected.
 Moreover, due to the existence of nonuniform convex constraint sets,
 we need take into account the nonlinearity of the consensus and optimal convergence caused by nonuniform convex constraint sets, which renders the analysis of this case to be very complicated. In particular when the nonuniform constraint sets are unbounded, {the gradients of the local objective functions might tend to be unbounded when their independent variables tend to infinity, which makes the existing approaches invalid}, e.g., \cite{angelia}, \cite{linren3}, where the nonuniform constraint sets and the subgradients were both bounded.
To solve the optimization problem, we perform the analysis in three steps. The first step is to make full use of the convexity of the objective functions and show that our algorithm ensures that all agents remain in a bounded region for all the time. The second step is to analyze the agent dynamics at some key times and show the consensus convergence of all agents. The third step is to estimate the consensus errors and use the estimation of the distance from the agents to the convex constraint sets and the convexity of the objective functions to show the optimal convergence of the optimization problem.
After that, we extend the obtained results to discrete-time multi-agent systems, and then study the case where each agent remains in its corresponding convex constraint set. Due to the discretization of the system dynamics, the agents might deviate from the minimum of the objective functions. Such a problem also exists in the centralized optimization system. To deal with it, a switching mechanism is introduced in the algorithms based on each agent's own information under which all agents remain in a bounded region even when the convex constraint sets are unbounded. It is shown that the distributed optimization problems with nonuniform convex constraint sets and nonuniform stepsizes can be solved for the discrete-time multi-agent systems.

This paper takes nonuniform unbounded convex constraint sets, nonuniform stepsizes, general differentiable convex objective functions, general switching graphs and the discretization of the algorithms into account simultaneously for the distributed optimization problems. The nonlinearities caused by these factors are different and the coexistence of these nonlinearities would further result in more complicated nonlinearities. Existing works only addressed a fraction of these factors due to the limitations of the algorithms and the analytical approaches. For example, the algorithms in \cite{cotes,Cortes3} cannot be directly applied to the case of convex constraint sets due to the adoption of the integrator operator. {The analytical approaches in \cite{linren4,xie6} cannot be directly applied in this paper, because the nonsmooth sign functions are used in \cite{linren4} to account for inconsistency in gradients while some intermediate variables need be transimitted besides the agent states in \cite{xie6}. Neither feature is valid in the current paper as continuous functions are used and no intermediate variables are transmitted. Moreover, in \cite{linren4,xie6}, the communication graphs are assumed to be strongly connected at all time or the constraint sets are assumed to be identical, which makes the analytical approaches in \cite{linren4,xie6} unable to be directly applied in this paper as well.}

Notation: $\mathbb{R}^m$  denotes the set of all $m$ dimensional real column
vectors; $\mathbb{R}^{m\times n}$  denotes the set of all $m\times
n$ real matrices; $\mathcal{I}$ denotes the index set
$\{1,\ldots,n\}$; $\textbf{1}$
represents
a column vector of all ones with a compatible dimension; $s_i$ denotes the $i$th component of the vector $s$; $A_{ij}$ denotes the $(i,j)$th entry of the matrix $A$;
$s^T$ and $A^T$ denote, respectively, the transpose of the vector $s$ and the matrix
$A$; $||s||$ denotes the Euclidean norm of the vector $s$; {$\nabla
f(s)$ denotes the gradient of the function $f(s)$ at $s$;} $\mathrm{diag}\{A_1,\cdots,A_n\}$ denotes a block diagonal matrices with its diagonal blocks equal to the matrices $A_i(k)$; 
  {the symbol $/$ denotes the division sign}; and {$P_{X}(s)$} denotes the projection of the vector $s$ onto the closed convex set $X$, i.e., {$P_{X}(s)=\mathrm{arg} \min\limits_{\bar{s}\in
X}\|s-\bar{s}\|$}.





\section{Preliminaries and problem formulation}
\subsection{Preliminaries}
Let $\mathcal{G}(\mathcal{I},\mathcal{E},\mathcal{A})$ be a
directed communication graph of agents, where $\mathcal{E}\subseteq\mathcal{I}\times
\mathcal{I}$ is the set of edges, and $\mathcal{A}=[a_{ij}]\in
\mathbb{R}^{n\times n}$ is the weighted adjacency matrix. An edge
$(i,j) \in \mathcal{E}$ denotes that agent $j$ can obtain
information from agent $i$. The weighted adjacency matrix $\mathcal{A}$ is
defined as $\eta\leq a_{ij}\leq \bar{\eta}$ for two constants $\bar{\eta}>\eta>0$ if $(j,i) \in
\mathcal{E}$ and $a_{ij}=0$ otherwise. It is assumed by default that $a_{ii}=0$, i.e., $(i,i)\notin \mathcal{E}$. The Laplacian of the graph $\mathcal{G}$, denoted by $L$, is defined as $L_{ii}=\sum_{j=1,j\neq i}^na_{ij}$ and $L_{ij}=-a_{ij}$ for all $i\neq j$. The graph $\mathcal{G}$ is undirected if $a_{ij}=a_{ji}$ for all $i,j$, and it is balanced if $\sum_{j=1}^na_{ij}=\sum_{j=1}^na_{ji}$ for all $i$. The set of neighbors of agent
$i$ is denoted by $N_i=\{j\in \mathcal{I}\mid (j,i)\in \mathcal{E}\}.$ A
path is a sequence of edges of the form
$(i_1,i_2),(i_2,i_3),\cdots$, where $i_j \in \mathcal{I}$. The graph $\mathcal{G}$ is strongly connected, if there
is a path from every agent to every other agent, and the graph $\mathcal{G}$ is connected, if it is undirected and strongly connected \cite{s10}.

\begin{lemma}\label{le13}{\rm\cite{s10} If the graph $\mathcal{G}$ is strongly connected, the Laplacian $L$ has one zero eigenvalue associated with eigenvalue vector $\mathbf{1}$ and  all its rest $n-1$ eigenvalues have positive real parts.
Further, if the graph $\mathcal{G}$ is undirected and connected, all the $n-1$ nonzero eigenvalues are positive.
}\end{lemma}
\begin{lemma}\label{lemma35}{\rm\cite{boyd} Let $f_0(\chi): \mathbb{R}^r\rightarrow\mathbb{R}$ be a differentiable convex function. $f_0(\chi)$ is minimized if and only if $\nabla f_0(\chi)=0$.}\end{lemma}

\begin{lemma}\label{le1u}{\rm \cite{Facchinei} Suppose that {$Y\neq\emptyset$} is a closed convex set in {$\mathbb{R}^r$}. The following statements hold.\\
(1)  For any $y\in \mathbb{R}^r$, $\|y-P_Y(y)\|$ is continuous with respect to $y$ and $\nabla \frac{1}{2}\|y-P_Y(y)\|^2=y-P_Y(y)$;\\
(2)  For any {$y,z\in\mathbb{R}^r$} and
all {$\xi\in Y$}, $[y-P_Y(y)]^T(y-\xi)\geq0$,
{$\|P_Y(y)-\xi\|^2\leq\|y-\xi\|^2-\|P_Y(y)-y\|^2$}
and {$\|P_Y(y)-P_Y(z)\|\leq\|y-z\|.$}}\end{lemma}

\subsection{Problem formulation}
Consider a multi-agent system consisting of $n$ agents. Each agent
is regarded as a node in a directed graph $\mathcal{G}${\footnote{In the following, $\mathcal{G}(t)$, $\mathcal{G}(kT)$, $a_{ij}(t)$, $a_{ij}(kT)$, $N_i(t)$ and $N_i(kT)$ will be used to denote the graph, the edge weight and the agent neighbor set at time $t$ or $k$ as defined in Sec. II.A.}}, and each agent can interact with only its local neighbors in $\mathcal{G}$.
Our objective is to design {{algorithms using only local interaction and information such that all agents cooperatively find the optimal state that solves
the optimization problem}}
\begin{eqnarray*}\label{gels1}\begin{array}{lll}\mathrm{minimize}~~\sum_{i=1}^nf_i(s)\\
\mathrm{subject~to}~~s\in H=\bigcap_{i=1}^n{H}_i\subseteq
\mathbb{R}^r,\end{array}\end{eqnarray*}
where $f_i(s): \mathbb{R}^r\rightarrow\mathbb{R}$ denotes the differentiable convex local objective
function of agent  $i$, and $H_i$ denotes the closed convex constraint set of $f_i(s)$. Clearly, $\sum_{i=1}^nf_i(s)$ is also a differentiable convex function. It is assumed that $f_i(s)$ and $H_i$ are known only to  agent
$i$. {{The problem described above is equivalent to the problem that all agents reach a consensus while minimizing the team objective function $\sum_{i=1}^nf_i(x_i)$, i.e.,
\begin{eqnarray}\label{gel1}\begin{array}{lll}\mathrm{minimize}~~\sum_{i=1}^nf_i(x_i)\\
\mathrm{subject~to}~~x_i=x_j\in
H=\bigcap_{i=1}^n{H}_i\subseteq
\mathbb{R}^r.\end{array}\end{eqnarray}}}
\hspace{-0.1cm}{{In this paper, our analysis is for the general $m$ case. When no confusion arises, the equations or formula are written in the form of $m=1$ for notational simplicity.}}





\section{Distributed continuous-time optimization}
{In this section, we discuss the distributed optimization problem for continuous-time multi-agent systems. The problem has applications in motion coordination
of multi-agent systems, where multiple physical vehicles rendezvous or form a formation centered
at a team optimal location. }
Suppose
that the agents satisfy the continuous-time dynamics
\begin{equation} \label{eq:single-int}
\dot{x}_i(t)=u_i(t), \quad i\in \mathcal{I},
\end{equation}
where $x_i\in \mathbb{R}^r$ is the state of agent $i$, and $u_i \in
\mathbb{R}^r$ is the control input of agent $i$.

\subsection{Assumptions and some necessary lemmas}
 {Let $\mathcal{X}$ denote the optimal set of the constrained optimization problem (\ref{gel1}). Before the main assumptions and the necessary lemmas, we need further define the sets $X_i\triangleq{\Big\{}x{\Big |}\nabla f_i(x)=0{\Big\}}$ for all $i\in\mathcal{I}$ and $X\triangleq{\Big\{}x{\Big |}\sum_{i=1}^n\nabla f_i(x)=0{\Big\}}$ for later usage. From Lemma \ref{lemma35}, $X_i$ and $X$ are convex and denote, respectively, the optimal sets of $f_i(x)$ and the team objective function $\sum_{i=1}^nf_i(x)$ for $x\in \mathbb{R}^r$. Note that in general $X$ is different from $\mathcal{X}$ but $X=\mathcal{X}$ when $H_i =\mathrm{R}^r$ for all $i$.}

\begin{assumption}\label{ass2}{\rm \cite{linren4} Each set $X_i$, $i\in\mathcal{I}$, is nonempty and bounded.
%
}\end{assumption}

{In Assumption \ref{ass2}, we only make an assumption on each $f_i(x)$ rather than the team objective function $\sum_{i=1}^nf_i(x)$ because $\sum_{i=1}^nf_i(x)$ is global information for all agents and cannot be used by each agent in a distributed way.}

\begin{assumption}\label{ass17}{\rm $H\neq\emptyset$.}\end{assumption}

{In Assumption \ref{ass17}, we do not require $H_i$ to be bounded and hence $H$ might be unbounded.}

\begin{lemma}\label{lemmallin}{\rm\cite{linren4} Under Assumption \ref{ass2},  the following two statements hold:\\
\noindent(1) $\lim_{\|y\|\rightarrow+\infty}f_i(y)=+\infty$ for all $i$ and accordingly $\lim_{\|y\|\rightarrow+\infty}\sum_{i=1}^nf_i(y)=+\infty$.

\noindent(2) All $X_i$, $i\in\mathcal{I}$, and $X$ are nonempty closed bounded convex sets.}\end{lemma}

\begin{lemma}\label{lemma109a1}{\rm \cite{linren4} Under Assumptions \ref{ass2} and  \ref{ass17}, $\mathcal{X}$ is a nonempty closed bounded convex set.}\end{lemma}

{Actually, Lemma \ref{lemma109a1} shows the existence and boundedness of the optimal set of the constrained optimization problem (\ref{gel1}), $\mathcal{X}$, under Assumptions \ref{ass2} and  \ref{ass17}.}

 \begin{assumption}\label{ass24}{\rm \cite{linren4} The length of the time interval between two contiguous switching times is no smaller than a given constant, denoted by $d_w$. }\end{assumption}

 {Under Assumption \ref{ass24}, the switching of the graph $\mathcal{G}(t)$ cannot be arbitrarily fast, which prevents the system from exhibiting the Zeno behavior.}

\begin{assumption}\label{ass16}{\rm There exists an infinite sequence of swiching
times of the graph $\mathcal{G}(t)$, {$t_0,t_1,t_2,\cdots,$} such that {$t_0=0,$ $0<t_{j+1}-t_j\leq M$} and the union of all the graphs during each interval $[t_j,t_{j+1})$
is strongly connected for some constant $M>d_w$ and all nonnegative integers $j$.}\end{assumption}

{Assumption \ref{ass16} ensures that all agents can communicate with each other persistently.}
 {Suppose that the graph $\mathcal{G}(t)$ is balanced for all $t$. From \cite{cortes00}, by rearranging the agent indices, $\frac{L(t)^T+L(t)}{2}$ can be written as $\frac{L(t)^T+L(t)}{2}=\mathrm{diag}\{\frac{L^1(t)^T+L^1(t)}{2},\cdots,\frac{L^h(t)^T+L^h(t)}{2}\}$, where each $\frac{L^i(t)^T+L^i(t)}{2}$ corresponds to a strongly connected component of the agents. From Lemma \ref{le13}, each eigenvalue of $\frac{L^i(t)^T+L^i(t)}{2}$ is nonnegative and hence all eigenvalues of $\frac{L(t)^T+L(t)}{2}$ are nonnegative.}

  Before the main results, we first present some necessary lemmas that will be used in the analysis of the main results. Specifically, Lemma \ref{lemma51} shows a  radial growth property of the derivatives of $f_i(x)$,  
 Lemma \ref{lemma9} shows the consensus convergence property of the stepsizes, Lemma \ref{lemma7} shows a dependency relationship between the distances from one given point to the convex sets and their intersection, and Lemma \ref{lemma10k} shows the boundedness of the gradients in a bounded region. For clarity, the proofs of Lemmas \ref{lemma9}-\ref{lemma10k} together with Theorems \ref{lemma1}-\ref{theorem2} are provided in the Appendix.

\begin{lemma}\label{lemma51}{\rm \cite{linren4} Let $f(s): \Xi\mapsto \mathbb{R}$ be a differentiable convex function and $Y$ be its minimum set in $\Xi$, where $\Xi\subseteq\mathbb{R}^r$ is a closed convex set. Suppose that $Y\subseteq \Xi$ is closed and bounded. For any $z=\lambda P_{Y}(y)+(1-\lambda)y$ with $\lambda\in(0,1)$, $0<\nabla f(z)^T\frac{y-P_{Y}(y)}{\|y-P_{Y}(y)\|}\leq \nabla f(y)^T\frac{y-P_{Y}(y)}{\|y-P_{Y}(y)\|}$ for any $y\in \Xi-Y$. }\end{lemma}


\begin{lemma}\label{lemma9}{\rm For the system given by $\dot{q}_i(t)=\mathrm{arctan} (e^{\|x_i(t)\|})$ with $q_i(0)>0$, if $\lim_{t\rightarrow+\infty}[x_i(t)-x^*(t)]=0$ for all $i$, $\lim_{t\rightarrow+\infty}\frac{q_i(t)}{q_j(t)}=1$ for all $i,j$, where $x^*(t)\triangleq\frac{1}{n}\sum_{i=1}^nx_i(t)$.}\end{lemma}





\begin{lemma}\label{lemma7}{\rm Let $y(t)\in E$ for all $t$, where $E$ is a bounded set. Under Assumption \ref{ass17}, if $\lim_{t\rightarrow+\infty}\|y(t)-P_{H_i}(y(t))\|=0$ for all $i$, then $\lim_{t\rightarrow+\infty}\|y(t)-P_{H}(y(t))\|=0$. }\end{lemma}

\begin{lemma}\label{lemma10k}{\rm Let $Y$ be a closed bounded convex set. Then, $\|\nabla f_i(x)\|<\varrho$ for all $i$, all $x\in Y$ and some constant $\varrho>0$.}\end{lemma}

\subsection{Algorithm and convergence analysis}

In this subsection, we design a continuous distributed optimization algorithms with nonuniform stepsizes. The algorithm is given by
   \begin{eqnarray}\label{eq5330f}\begin{array}{lll}
\dot{q}_i(t)&=&\mathrm{arctan} (e^{\|x_i(t)\|}), q_i(0)>0,\\
u_i(t)&=&\sum_{j\in
N_i(t)}a_{ij}(t){\big(}x_j(t)-x_i(t){\big)}\\&-&[x_i(t)-P_{H_i}(x_i(t))]-\frac{\nabla
f_i(x_i(t))}{\sqrt{q_i(t)}}\end{array}
\end{eqnarray}
for all $i$. The stepsize of the gradient,  $1/\sqrt{q_i(t)}$, is used to make the term $\frac{\nabla
f_i(x_i(t))}{\sqrt{q_i(t)}}$ tend to zero as $t\rightarrow+\infty$. The role of the term $-\frac{\nabla
f_i(x_i(t))}{\sqrt{q_i(t)}}$ is to make all agents converge to the optimal set of the team objective function and the role of $-[x_i(t)-P_{H_i}(x_i(t))]$ is to make each agent converge to the convex set $H_i$.

{It should be noted that the construction of the stepsize $1/\sqrt{q_i(t)}$ is only based on the $i$th agent's own states and it does not use the Lipschitz constant or the convexity constant as in the existing works, e.g., \cite{cotes} and \cite{hong4}.}
The stepsize of each agent $1/\sqrt{q_i(t)}$ is state-dependent and can be nonuniform for all agents. The existing works often assume the stepsizes of all agents to be predesigned and consistent with each other at any time, which exerts a heavy burden on sensing and communication costs of the entire system.

\begin{remark}{\rm
In algorithm (\ref{eq5330f}), the role of the inverse tangent functions and the exponential functions is to ensure $\dot{q}_i(t)$ to be upper and lower bounded. In fact, some other more general functions, e.g., saturation function, can be employed to play the same role. Moreover, the stepsizes used in algorithm (\ref{eq5330f})  are in a special form, and it can also be extended to other functions. For easy readability, we do not give the general form of the functions $q_i(t)$ and the stepsizes.
}\end{remark}

Let $$x(t)=[x_1(t)^T,\cdots,x_n^T(t)]^T,$$ $$Q(t)=\mathrm{diag}\{1/\sqrt{q_1(t)},\cdots,1/\sqrt{q_n(t)}\},$$ $$\nabla
f(x(t))=[\nabla
f_1(x_1(t))^T,\cdots,\nabla
f_n(x_n(t))^T]^T$$
and $Pz(x(t))=[x_1(t)-P_{H_1}(x_1(t)),\cdots,x_n(t)-P_{H_n}(x_n(t))]^T$. Then the system (\ref{eq:single-int}) with (\ref{eq5330f}) can be written as
\begin{eqnarray}\label{e901}\begin{array}{lll}\dot{x}(t)=-L(t){x}(t)-Pz(x(t))-Q(t)\nabla
f(x(t)).\end{array}\end{eqnarray}

 For convenience of discussion, under Assumption \ref{ass16}, let $t_{j_0}<t_{j_1}<\cdots<t_{j_{m_j}}$ with $t_{j_0}=t_j$ and $t_{j_{m_j}}=t_{j+1}$ denote all the switching times in the interval $[t_j,t_{j+1})$.
In the following, we show the effectiveness of this algorithm to solve the optimization problem (\ref{gel1}).
\begin{theorem}\label{lemma1}{\rm Suppose
that the graph $\mathcal{G}(t)$ is balanced for all $t$ and Assumptions \ref{ass2}, \ref{ass17}, \ref{ass24} and \ref{ass16} hold.
 For arbitrary initial conditions $x_i(0)\in \mathbb{R}^r$, using algorithm (\ref{eq5330f}) for system (\ref{eq:single-int}),
 the following statements hold.
\begin{itemize}
\item [(1)] All $x_i(t)$ remain in a bounded region for all $i$ and all $t$.
\item [(2)] $\lim_{t\rightarrow+\infty}\|x^*(t)-P_{H}(x^*(t))\|=\lim_{t\rightarrow+\infty}[x_i(t)-x^*(t)]=0$ for all $i$, where $x^*(t)$ has been defined in Lemma \ref{lemma9}.
\end{itemize}
}\end{theorem}



\begin{theorem}\label{theorem1}{\rm Suppose
that the graph $\mathcal{G}(t)$ is balanced for all $t$ and Assumptions  \ref{ass2}, \ref{ass17}, \ref{ass24} and \ref{ass16} hold.
 For arbitrary initial conditions $x_i(0)\in \mathbb{R}^r$, using algorithm (\ref{eq5330f}) for system (\ref{eq:single-int}), all agents reach a
consensus, i.e., $\lim_{t\rightarrow+\infty}[x_i(t)-x^*(t)]=0$ for all $i$, and minimize the
 team objective function (\ref{gel1}) as $t\rightarrow+\infty$.
}
\end{theorem}

\begin{remark}{\rm In the existing works, the distributed optimization problems were considered often under the assumption that the stepsizes of the gradients are uniform and explicitly time-dependent for all agents. That is, the stepsizes should be consistent with each other at any time.
In algorithm (\ref{eq5330f}), we do not make such an assumption and the stepsizes of the gradients only depend on the agent states. The stepsizes need not have the same value at any time and instead they are usually nonuniform, which greatly relaxes the synchronization requirement on the system.
}\end{remark}


\section{Distributed discrete-time optimization}
In this section, we discuss the distributed optimization problem for discrete-time multi-agent systems.
Suppose
that the agents satisfy the discrete-time dynamics
\begin{equation} \label{eq:single-int1}
{x}_i((k+1)T)=u_i(kT), \quad i\in \mathcal{I},
\end{equation}
where $x_i\in \mathbb{R}^r$ is the state of agent $i$, $u_i \in
\mathbb{R}^r$ is the control input of agent $i$, and $T>0$ is the sample time. In the following, we use ``$(k)$" instead of ``$(kT)$" when no confusion arises.

\subsection{Assumptions and some necessary lemmas}
In Sec. II, when we define the weighted adjacency matrix $\mathcal{A}$, we assume by default that $a_{ii}=0$. In this section, for discussion of discrete-time multi-agent systems, we need to redefine $a_{ii}$ and make an assumption about the weighted adjacency matrix $\mathcal{A}$ as shown in the following assumption.
\begin{assumption}\label{ass982}{\rm \cite{linren3} For all {$i,j\in \mathcal{I}$}, $a_{ii}(k)\geq \eta$, $a_{ij}(k)\geq \eta$ for some constant $0<\eta\leq1$ and each nonzero $a_{ij}(k)$, {$\sum_{i=1}^na_{ij}(k)=1$} and
{$\sum_{j=1}^n a_{ji}(k)=1$}. }\end{assumption}

{Under Assumption \ref{ass982}, the adjacency matrix of the graph $\mathcal{G}(k)$ is doubly stochastic and its diagonal entries are nonzero. Assumption \ref{ass982}
is used to generate convex combinations of the agents' states such that the influence of each agent's state is equal in the final consensus value in the distributed optimization algorithms shown later.}

\begin{assumption}\label{ass16o}{\rm\cite{linren3}
There exists an infinite sequence of switching
times of the graph $\mathcal{G}(k)$, {$k_0,k_1,k_2,\cdots,$} such that $k_0=0$, $0<k_{j+1}-k_j\leq M$ and the union of all the graphs during each interval $[k_j,k_{j+1})$
is strongly connected for some positive integer $M$ and all nonnegative integers $j$.}\end{assumption}

Similar to Assumption \ref{ass16}, {Assumption \ref{ass16o} ensures that all agents can communicate with each other persistently.}



\subsection{Distributed optimization with nonuniform stepsizes}
In this subsection, we design a discrete-time distributed optimization algorithm with nonuniform stepsizes. The algorithm is given by \begin{eqnarray}\label{eq53371}\begin{array}{lll}
&{q}_i(k+1)={q}_i(k)+\mathrm{arctan} (e^{\|x_i(k)\|})T, q_i(0)>0,\\
&u_i(k)=w_i(k)(1-\gamma_i)+P_{H_i}[w_i(k)]\gamma_i,\\
&w_i(k)=v_i(k)-gr_i(k)T,\\
&v_i(k)=\sum_{j\in
N_i(k)\cup \{i\}}a_{ij}(k)x_j(k),\\
&gr_i(k)=\left\{\begin{array}{lll}0, \mbox{~if~} \sqrt{q_i(k)}\leq\|\nabla
f_i(v_i(k))\|^2,\\
\frac{\nabla
f_i(v_i(k))}{\sqrt{q_i(k)}}, \mbox{~otherwise,}\end{array}\right.\end{array}
\end{eqnarray}
where  $0<\gamma_i\leq1$ is a constant for each $i$.

{Due to the discretization of the system dynamics, the agents might deviate from the minimum of the team objective function. Such a problem also exists in the centralized optimization system. To deal with it, a switching mechanism is introduced in (\ref{eq53371}) based on each agent' own information, under which the gradient term would not be too large to result in the divergence of the system. This will be shown in the proof of Theorem \ref{theorem14e}.}

{For Algorithm (\ref{eq53371}), it can be calculated simply in four steps: (a) $v_i(k)$, $q_i(k+1)$ and $\nabla f_i(v_i(k))$; (b) $gr_i(k)$ based on the switching mechanism; (c) $w_i(k)$; and (d) $P_{H_i}[w_i(k)]$ and $u_i(k)$. Though the algorithm computation looks a bit complex due to the existence of the switching mechanism, the algorithm does not require intermediate variables to be transmitted and it is a fully distributed algorithm.}

\begin{theorem}\label{theorem14e}{\rm Suppose
that  Assumptions \ref{ass2}, \ref{ass17}, \ref{ass982} and \ref{ass16o} hold.
For arbitrary initial conditions $x_i(0)\in \mathbb{R}^r$, using algorithm (\ref{eq53371}) for system (\ref{eq:single-int1}), if $0<\gamma_i<1$ for all $i$, all agents reach a
consensus, i.e., $\lim_{t\rightarrow+\infty}[x_i(k)-x^*(k)]=0$ for all $i$, and minimize the
 team objective function (\ref{gel1}) as $k\rightarrow+\infty$.
}
\end{theorem}

In Theorems \ref{theorem1} and \ref{theorem14e}, it is not required that each agent remain in its corresponding convex constraint set $H_i$. In the following theorem, we show that the optimization problem (\ref{gel1}) can be solved when all agents remain in their corresponding convex constraint sets.

\begin{theorem}\label{theorem2}{\rm Suppose
that  Assumptions \ref{ass2}, \ref{ass17}, \ref{ass982} and \ref{ass16o} hold.
For arbitrary initial conditions $x_i(0)\in H_i$, using algorithm (\ref{eq53371}) for system (\ref{eq:single-int1}), if $\gamma_i=1$ for all $i$, all agents reach a
consensus, i.e., $\lim_{t\rightarrow+\infty}[x_i(k)-x^*(k)]=0$ for all $i$, and minimize the
 team objective function (\ref{gel1}) as $k\rightarrow+\infty$ while each agent remains in its corresponding constraint sets, i.e., $x_i(k)\in H_i$ for all $i$ and all $k$.
}
\end{theorem}


\begin{remark}{\rm Since the proposed algorithms are gradient based, the convergence rate of the algorithms is not very fast. This is common for the gradient-based distributed algorithms in the existing literature. In particular, the stepsizes (gradient gains) are nonuniform, which makes the convergence rate slower than that with uniform stepsizes. However, our algorithms are able to deal with the general case of nonuniform stepsizes without intermediate variables being transmitted. In the existing non-gradient-based works, some special assumptions are always made in order to ensure the optimal convergence. For example, in \cite{Wei}, the communication graphs are assumed to be strongly connected and the the convex objective functions are assumed to be strongly convex. In this paper, the communication graphs are only required to be jointly strongly connected, and the convex objective functions are only required to be differentiable (which can be easily extended to the nondifferentiable case by using subgradients). Future work could be directed towards improving the convergence rate of our algorithms. In particular, different dimension might yield different convergence rates. It is worth studying the effects of different dimensions on the convergence rate of the algorithms.}\end{remark}

\section{Simulations}
{Consider a multi-agent system with
8 agents in $\mathbb{R}^2$.  The communication graphs switch among the balanced subgraphs of the graph shown in Fig. \ref{fig:1}.
 Each edge weight is 0.5. The
sample time is $T=0.1$ $s$ for the discrete-time algorithm. The local objective
functions are adopted as $f_{1}(x_1)=\frac{1}{2}x_{11}^2+\frac{1}{2}x_{12}^2,$ $f_2(x_2)=\frac{1}{2}(x_{21}+1)^2+\frac{1}{2}x_{22}^2,$
$f_3(x_3)=\frac{1}{2}x_{31}^2+\frac{1}{2}(x_{32}+1)^2,$
$f_4(x_4)=\frac{1}{2}(x_{41}+1)^2+\frac{1}{2}(x_{42}+1)^2$, $f_5(x_5)=\frac{1}{4}x_{51}^4+\frac{1}{4}x_{52}^4,$
$f_6(x_6)=\frac{1}{4}(x_{61}+1)^4+\frac{1}{4}x_{62}^4,$
$f_7(x_7)=\frac{1}{4}x_{71}^4+\frac{1}{4}(x_{72}+1)^4$, and $f_8(x_8)=\frac{1}{4}(x_{81}+1)^4+\frac{1}{4}(x_{82}+1)^4$.
 where $x_{i1}$ and $x_{i2}$ denote the $1$st and $2$nd components of $x_i$. The constrained convex sets are adopted as $H_1=\{(x,y)^T\in \mathbb{R}^2 \mid \|s-[0,0]^T\|\leq 3\}$ for agents 1 and 5, $H_2=\{(x,y)^T\in \mathbb{R}^2 \mid \|x\leq 0.5, y\geq 1\}$ for agents 2 and 6, $H_3=\{(x,y)^T\in \mathbb{R}^2 \mid \|s-[0,3]^T\|\leq 3\}$ for agents 3 and 7 and $H_4=\{(x,y)^T\in \mathbb{R}^2 \mid x\geq -0.5,y\geq 1\}$ for agents 4 and 8.
The team objective function $\sum_{i=1}^8f_i(s)$ is minimized if and only if $s=[-0.5,1]^T$. The simulation
results for algorithm (\ref{eq5330f}) and algorithm (\ref{eq53371}) with $\gamma_i=1$ for all $i$ are shown in Figs. \ref{fig:2} and \ref{fig:3}. It is clear that all agents finally converge to the optimal point. In particular, for algorithm (\ref{eq53371}), all agents remain in their corresponding constraint sets. All the simulation results are
 consistent with our obtained theorems. 
\begin{figure}
\begin{center}
\hspace{0.2cm}\xymatrix{
    *++[o][F-]{1}\ar[r]\ar[d]
    & *++[o][F-]{2}\ar[l]\ar[r]
    &*++[o][F-]{3}\ar[l]\ar[r]
    & *++[o][F-]{4}\ar[d]\ar[l]
\\
    *++[o][F-]{8}\ar[u]\ar[r]
    & *++[o][F-]{7}\ar[l]\ar[r]
    & *++[o][F-]{6}\ar[r]\ar[l]
    & *++[o][F-]{5}\ar[l]\ar[u]
  }
\end{center}
\vspace{0.2cm}
\caption{One undirected graph.} \label{fig:1}
\end{figure}
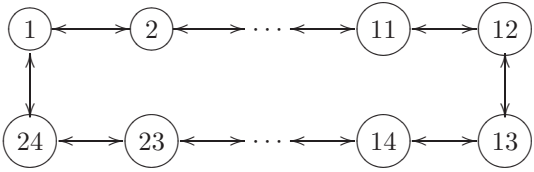

\begin{figure}
\centering
\hspace{-0.56cm}\includegraphics[width=3.7in]{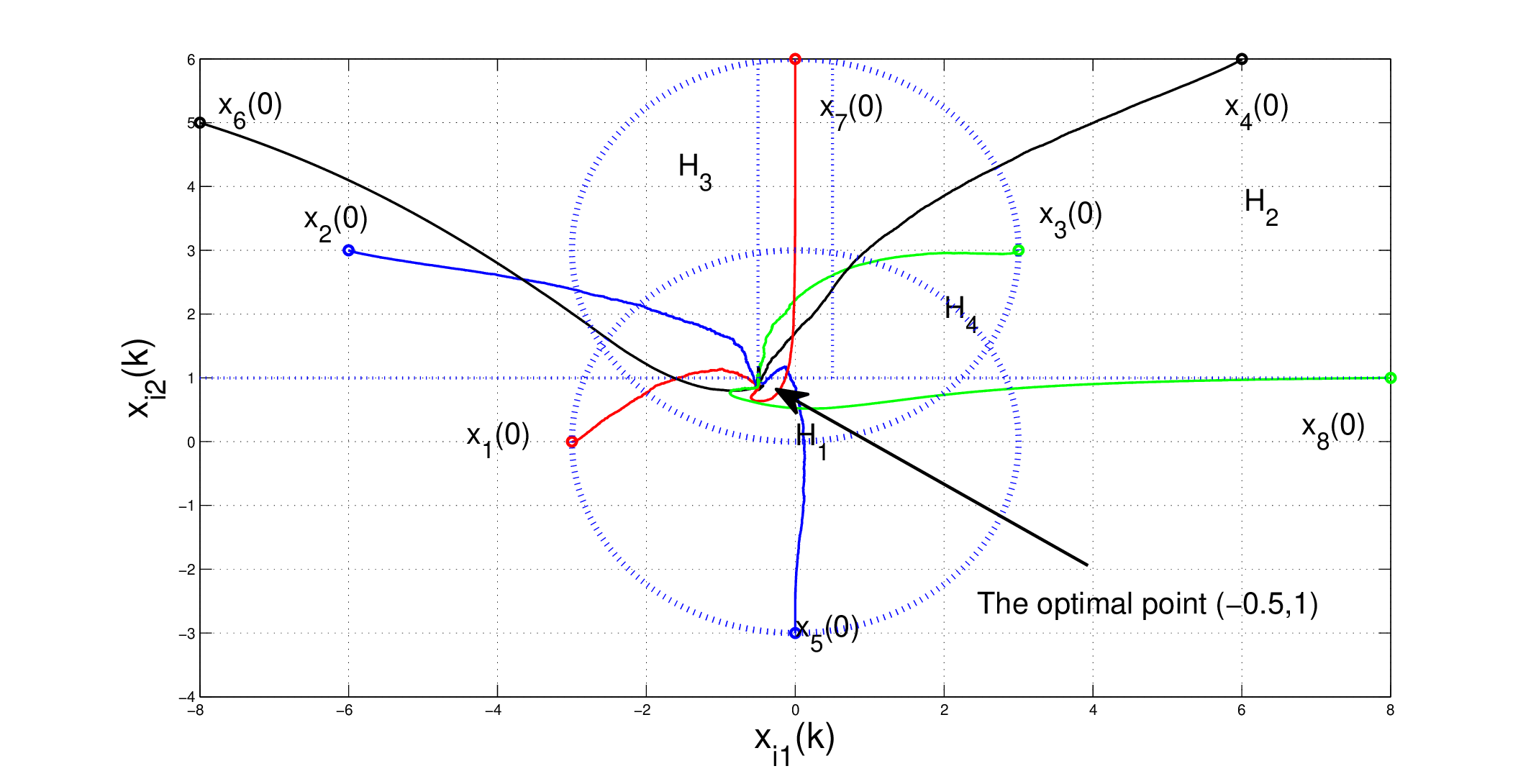} \\
\vspace{-0.2cm} \caption{State trajectories of all agents using (\ref{eq5330f}).}
\label{fig:2}
\end{figure}
\begin{figure}
\centering
\hspace{-0.56cm}\includegraphics[width=3.7in]{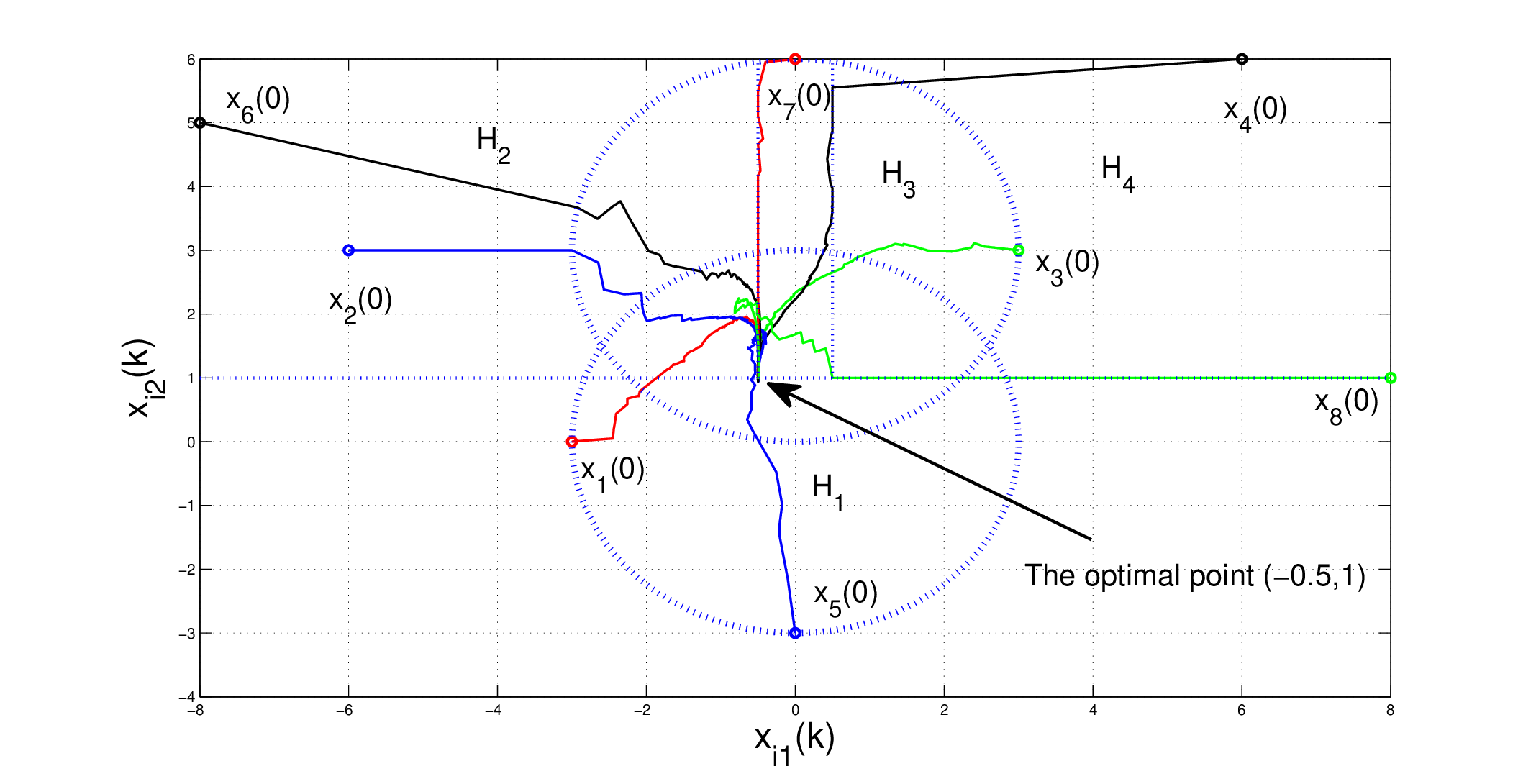} \\
\vspace{-0.2cm} \caption{State trajectories of all agents using (\ref{eq53371}).}
\label{fig:3}
\end{figure}}

\section{Conclusions}
In this paper, distributed continuous-time and discrete-time optimization problems with nonuniform stepsizes and nonuniform possibly unbounded convex constraint sets were studied
for general differentiable convex objective functions. One continuous-time algorithm and one discrete-time algorithm were introduced.
For these two algorithms, it has been shown that
the team objective function is minimized with nonuniform stepsizes and nonuniform convex constraint sets. In particular, for the discrete-time algorithm, it has been shown that the distributed optimization problem can be solved when each agent remain in its corresponding constraint set.

\appendix
\noindent{\textbf{Proof of Lemma \ref{lemma9}:}} Note that $\pi/4\leq\mathrm{arctan}(e^{\|x_i(t)\|})\leq \pi/2$ for all $i$ and all $t$. There exists a constant $T_0>0$ such that $\frac{t}{4}<q_i(t)<4t$ for all $t>T_0$ and all $i$.
Let $\Delta_i(t)=q_i(0)-q_1(0)+\int_{0}^t[\mathrm{arctan}(e^{\|x_i(s)\|})-\mathrm{arctan}(e^{\|x_1(s)\|})]\mathrm{d}s$. It is clear that $q_i(t)=q_1(t)(1+\Delta_i(t)/q_1(t))$ for all $i$. Since $\lim_{t\rightarrow+\infty}[x_i(t)-x^*(t)]=0$, from the continuity of the function $\mathrm{arctan}(e^{\|x_i(t)\|})$, there exists a constant $T_1>T_0$ for any ${\epsilon}_0>0$ such that
 $|\mathrm{arctan}(e^{\|x_i(t)\|})-\mathrm{arctan}(e^{\|x_1(t)\|})|<{\epsilon}_0$ for all $t>T_1$ and all $i$.
It is clear that $\lim_{t\rightarrow+\infty}\Delta_i(T_1)/q_1(t)=0$ and $\lim_{t\rightarrow+\infty}|\int_{T_1}^{t}[\mathrm{arctan}(e^{\|x_i(s)\|})-\mathrm{arctan}(e^{\|x_1(s)\|})]\mathrm{d}s/q_1(t)|$ $<4{\epsilon}_0$. Since ${\epsilon}_0$ can be arbitrarily chosen, we have
\begin{eqnarray*}\label{er1}&&\lim_{t\rightarrow+\infty}\Delta_i(t)/q_1(t)=\lim_{t\rightarrow+\infty}[\Delta_i(T_1)\\&+&\int_{T_1}^{t}[\mathrm{arctan}(e^{\|x_i(s)\|})-\mathrm{arctan}(e^{\|x_1(s)\|})]\mathrm{d}s]/q_1(t)=0\end{eqnarray*} and hence $\lim_{t\rightarrow+\infty}\frac{q_i(t)}{q_1(t)}=1$ for all $i$. Therefore, $\lim_{t\rightarrow+\infty}\frac{q_i(t)}{q_j(t)}=1$ for all $i,j$.\endproof

\noindent{\textbf{Proof of Lemma \ref{lemma7}:}}  Under Assumption \ref{ass17}, $H\neq \emptyset$. Since $E$ is a bounded set, there exists a closed bounded convex set $Y$ such that $E\subset Y$. From Lemma \ref{le1u}, $\|y-P_{H}(y)\|$ is a continuous function with respect to $y$. If this lemma does not hold, there must exist a sequence $\{t^k\}$ such that $\lim_{k\rightarrow+\infty}\|y(t^k)-P_{H}(y(t^k))\|=\psi$ for some constant $0<\psi\leq+\infty$. Let $H_{\psi}=\{y\mid \|y-P_{H}(y)\|=\psi\}$. Consider the set $H_{\psi}\cap Y$. If $\sum_{i=1}^n\|y-P_{H_i}(y)\|=0$ for some $y$, we have $y\in H_i$ for all $i$ and hence $y\in H$, i.e., $\|y-P_{H}(y)\|=0$.

Note that  $H_{\psi}\cap Y$ is a closed bounded set. The inequality $\|y-P_{H}(y)\|=\psi>0$ implies that
 the function $\sum_{i=1}^n\|y-P_{H_i}(y)\|$ on $H_{\psi}\cap Y$ is lower bounded by a positive constant, denoted by $\gamma$.
  If this is not true, there must exist a
 sequence $\{\hat{y}_h\in H_{\psi}\cap Y\}$ such that $\lim_{h\rightarrow+\infty}\sum_{i=1}^n\|\hat{y}_h-P_{H_i}(\hat{y}_h)\|=0$.
Note that the sequence $\{\hat{y}_h\in H_{\psi}\cap Y\}$ has a subsequence $\{\bar{y}_h\in H_{\psi}\cap Y\}$ such that $\lim_{h\rightarrow+\infty}\bar{y}_h=\bar{y}$ for a constant vector $\bar{y}$. Since $H_{\psi}\cap Y$ is a closed bounded set, $\bar{y}\in H_{\psi}\cap Y$, i.e., $\sum_{i=1}^n\|\bar{y}-P_{H_i}(\bar{y})\|=0$, and hence $\|\bar{y}-P_{H}(\bar{y})\|=0$. This contradicts with the definition of $H_{\psi}$.

 Since $y(t)\in Y$ for all $t$, $\lim_{k\rightarrow+\infty}\|y(t^k)-P_{H}(y(t^k))\|=\psi$ and $\lim_{t\rightarrow+\infty}\|y(t)-P_{H_i}(y(t))\|=0$ for all $i$, there exists a constant $N>0$ for any $\epsilon>0$ such that $\sum_{i=1}^n\|y(t^k)-P_{H_i}(y(t^k))\|-\gamma>-\epsilon$ and $\|y(t)-P_{H_i}(y(t))\|<\epsilon$ for all $t^k>N$ and $t>N$. Let $\epsilon<\frac{\gamma}{n+1}$. It follows that $\sum_{i=1}^n\|y(t^k)-P_{H_i}(y(t^k))\|>\gamma-\epsilon>n\epsilon>\sum_{i=1}^n\|y(t)-P_{H_i}(y(t))\|$ for all $t^k>N$ and $t>N$, which is a contradiction.
 \endproof

 \noindent{\textbf{Proof of Lemma \ref{lemma10k}:}} Since each function $f_i(x)$ is differentiable, $\nabla f_i(x)$ exists for any $x\in Y$. If $\|\nabla f_i(x)\|$ is unbounded for some $x\in Y$, there must exist a sequence $\{y^k\in Y,k=1,2,\cdots\}$ such that $\lim_{k\rightarrow+\infty}y^k=y^0$ for some $y^0\in Y$ and $\lim_{k\rightarrow+\infty}\|\nabla f_i(y^k)\|=+\infty$. Let $z^k=y^k+\nabla f_i(y^k)/\|\nabla f_i(y^k)\|$. It is clear that each $f_i(z^k)$ and each $f_i(y^k)$ are bounded. From the convexity of $f_i(y^k)$, it follows that
$\|\nabla f_i(y^k)\|=\nabla f_i(y^k)^T(z^k-y^k)\leq f_i(z^k)-f_i(y^k)$ and hence $\lim_{k\rightarrow+\infty}[f_i(z^k)-f_i(y^k)]\geq +\infty$ which is a contradiction.\endproof

\noindent{\textbf{Proof of Theorem \ref{lemma1}:}}  First, we prove that using (\ref{eq5330f}) for (\ref{eq:single-int}) all $x_i(t)$ remain in a bounded region and all $\|\nabla f_i(x_i(t))\|$ are bounded for all $i$ and all $t$. 
  Note that $\pi/4\leq\mathrm{arctan}(e^{\|x_i(t)\|})\leq \pi/2$ for all $t$ and all $i$.
There exists a constant $T_0>0$ such that $2\sqrt{t}>\sqrt{q_i(t)}>\frac{\sqrt{t}}{2}$ for all $i$ and all $t\geq T_0$. Under Assumption \ref{ass2}, from Lemma \ref{lemmallin}, $X$ and all $X_i$ are nonempty and bounded. Note from Lemma \ref{lemma109a1} that $\mathcal{X}$ is nonempty and bounded under Assumptions \ref{ass2} and \ref{ass17}. There is a closed bounded convex set $Y=\{y\mid \|y-P_{\mathcal{X}}(y)\|\leq P_1\}\subset \mathbb{R}^r$ for some constant $P_1>0$ such that $x_i(T_0)\in Y$,
$X\subset Y$, $X_i\subset Y$ and $\mathcal{X}\subset Y$ for all $i$. From the continuity of $f_i(x)$ and Lemma \ref{lemmallin}, let $Y$ be sufficiently large such that $f_i(x_i(t))-f_i(z)\geq 4\sum_{j=1}^n[f_j(z)-f_j(z_j)]$ for all $i$, all $z\in \mathcal{X}$, all $z_j\in X_j$ and all $x_i(t)\notin Y$.

Consider the Lyapunov function candidate
$V_0(t)=\frac{1}{2}\|x(t)-z \mathbf{1}\|^2$ for some $z\in \mathcal{X}$ and $t\geq T_0$. Calculating $\dot{V}_0(t)$, we have
\begin{eqnarray*}\begin{array}{lll}&\dot{V}_0(t)\\&={\big(}x(t)-z\mathbf{1}{\big)}^T[-L(t){x}(t)-Pz(x(t))-Q(t)\nabla
f(x(t))]
\end{array}\end{eqnarray*}for $t\geq T_0$.
Since the graph $\mathcal{G}(t)$ is balanced, then $\mathbf{1}^TL(t)=0$ and hence $z\mathbf{1}^TL(t){x}(t)=0$. Regarding $\frac{L(t)^T+L(t)}{2}$ as a Laplacian of an undirected graph, we have $-x(t)^TL(t){x}(t)=-x(t)^T\frac{L(t)^T+L(t)}{2}{x}(t)\leq 0$ from the previous analysis below Assumption \ref{ass16}. Since the function $\frac{1}{2}\|x_i(t)-P_{H_i}(x_i(t))\|^2$ is convex, it follows that \begin{eqnarray*}\begin{array}{lll}-(x_i(t)-z)^T(x_i(t)-P_{H_i}(x_i(t)))\\\leq \frac{1}{2}\|z-P_{H_i}(z)\|^2-\frac{1}{2}\|x_i(t)-P_{H_i}(x_i(t))\|^2\\=-\frac{1}{2}\|x_i(t)-P_{H_i}(x_i(t))\|^2\end{array}\end{eqnarray*} where the last equality holds since $z\in \mathcal{X}\subset H\subset H_i$. Thus, we have
\begin{eqnarray*}\begin{array}{lll}\dot{V}_0(t)\leq-\frac{1}{2}\sum_{i=1}^n\|x_i(t)-P_{H_i}(x_i(t))\|^2\\-{\big(}x(t)-z\mathbf{1}{\big)}^TQ(t)\nabla
f(x(t))
\end{array}\end{eqnarray*}for $t\geq T_0$.
From the convexity of the function $f_i(x_i(t))$, we have \begin{eqnarray*}\begin{array}{lll}\nabla f_i(x_i(t))^T(z-x_i(t))\leq f_i(z)-f_i(x_i(t))\leq f_i(z)-f_i(z_i).\end{array}\end{eqnarray*}It follows that
\begin{eqnarray*}\begin{array}{lll}-{\big(}x(t)-z\mathbf{1}{\big)}^TQ(t)\nabla
f(x(t))\\\leq\sum_{i=1}^n\frac{1}{\sqrt{q_i(t)}}[f_i(z)-f_i(x_i(t))]. \end{array}\end{eqnarray*} Note that $f_i(x_i(t))-f_i(z_i)\geq0$ since $z_i\in X_i$ for all $i$. If $x_{i_0}(t)\notin Y$ for some $i_0$, we have $f_{i_0}(x_{i_0}(t))-f_{i_0}(z)\geq 4\sum_{j=1,j\neq i_0}^n[f_j(z)-f_j(z_j)]$ for all $z_j\in X_j$ and hence \begin{eqnarray*}\begin{array}{lll}-{\big(}x(t)-z\mathbf{1}{\big)}^TQ(t)\nabla
f(x(t))\\\leq -\frac{1}{\sqrt{q_{i_0}(t)}}[f_{i_0}(x_{i_0}(t))-f_{i_0}(z)]\\+\sum_{j=1,j\neq i_0}^n\frac{1}{\sqrt{q_j(t)}}[f_j(z)-f_j(z_j)]\\
\leq -\frac{1}{2\sqrt{t}}[f_{i_0}(x_{i_0}(t))-f_{i_0}(z)]\\+\sum_{j=1,j\neq i_0}^n\frac{2}{\sqrt{t}}[f_j(z)-f_j(z_j)]
\leq 0.\end{array}\end{eqnarray*}
 As a result, we have
$\dot{V}_0(t)\leq 0$ for $t\geq T_0$ if there exists an agent $i_0$ such that $x_{i_0}(t)\notin Y$ for some $i_0$. This implies that all $x_i(t)$ remain in a bounded region for all $i$ and all $t$. Thus, $\|x_i(t)-x^*(t)\|$ is bounded. {Note that $\|x_i(t)-P_{H}(x^*(t))\|\leq \|x^*(t)-P_{H}(x^*(t))\|+\|x_i(t)-x^*(t)\|\leq \|x_i(t)-x^*(t)\|+\|x^*(t)-s\|$ where $s\in H$ is a constant vector. Hence, $\|x_i(t)-P_{H}(x^*(t))\|$ is bounded.}
 Further, from Lemma \ref{lemma10k}, it follows that $\|\nabla f_i(x_i(t))\|$ is bounded for all $i$ and all $t$.  Since $2\sqrt{t}>\sqrt{q_i(t)}>\frac{\sqrt{t}}{2}$ for all $i$ and all $t\geq T_0$, we have  \begin{eqnarray*}\begin{array}{lll}\|x_i(t)-P_{H}(x^*(t))\|\|\nabla f_i(x_i(t))\|/\sqrt{q_i(t)}<\frac{\rho}{n\sqrt{t}}\end{array}\end{eqnarray*} for all $i$, all $t\geq T_0$ and some constant $\rho>0$. Moreover, since the graph $\mathcal{G}(t)$ is balanced, \begin{eqnarray*}\begin{array}{lll}\dot{x}^*(t)=-\frac{1}{n}\sum_{i=1}^n[x_i(t)-P_{H_i}(x_i(t))+\frac{\nabla
f_i(x_i(t))}{\sqrt{q_i(t)}}].\end{array}\end{eqnarray*}

Consider the Lyapunov function candidate
\begin{eqnarray}\label{lyap1}V_1(t)=\sum_{i=1}^n\|x_i(t)-x^*(t)\|^2+n\|x^*(t)-P_{H}(x^*(t))\|^2\end{eqnarray} for all $t\geq T_0$.
Calculating $\dot{V}_1(t)$, we have
\begin{eqnarray}\label{e121}\begin{array}{lll}&\dot{V}_1(t)\\=&2\sum_{i=1}^n[x_i(t)-x^*(t)]^T(\dot{x}_i(t)-\dot{x}^*(t))\\
+&2n[x^*(t)-P_{H}(x^*(t)]^T\dot{x}^*(t)\\
=&-2x^T(t)L(t)x(t)-2\sum_{i=1}^n[x_i(t)-x^*(t)]^T\\\times&[(x_i(t)-P_{H_i}(x_i(t))+\frac{\nabla
f_i(x_i(t))}{\sqrt{q_i(t)}}]\\
-&2\sum_{i=1}^n[x^*(t)-x_i(t)+x_i(t)-P_{H}(x^*(t))]^T\\\times&[x_i(t)-P_{H_i}(x_i(t))+\frac{\nabla
f_i(x_i(t))}{\sqrt{q_i(t)}}]\\
=&-2x^T(t)L(t)x(t)-2\sum_{i=1}^n[x_i(t)-P_{H}(x^*(t))]^T\\\times&[x_i(t)-P_{H_i}(x_i(t))+\frac{\nabla
f_i(x_i(t))}{\sqrt{q_i(t)}}]\\
\leq&-\phi(t)+\frac{2\rho}{\sqrt{t}}
\end{array}
\end{eqnarray}
for all $t\geq T_0$, where \begin{eqnarray*}\begin{array}{lll}\phi(t)\triangleq2x(t)^TL(t){x}(t)+\sum_{i=1}^n\|x_i(t)-P_{H_i}(x_i(t))\|^2, \end{array}\end{eqnarray*} $\sum_{i=1}^n[x_i(t)-x^*(t)]^T\dot{x}^*(t)=0$, and the last inequality has used the convexity of $\|x_i(t)-P_{H_i}(x_i(t))\|^2$ and the fact that $P_{H}(x^*(t))\in H_i$ and $\|P_{H}(x^*(t))-P_{H_i}(P_{H}(x^*(t)))\|^2=\|P_{H}(x^*(t))-P_{H}(x^*(t))\|^2=0$.

To proceed with our proof, we prove that $V_1(t)\leq \delta$ for any constant $\delta>0$ when $t$ is sufficiently large.
Since all $x_i(t)$ remain in a bounded region for all $i$ and all $t$, each $\|\dot{x}_i(t)\|$ is bounded for all $i,t$.
There exists a constant $h_0>0$ such that $\|\dot{\phi}(t)\|<h_0$ for all $t\geq T_0$.
 Consider $V_1(t)$ for $t\in [t_k,t_{k+1})$ and $t_k\geq T_0$ where $t_k$ has been defined in Assumption \ref{ass16}. Suppose that $h_1(t_k)=\max_{s\in [{t_k},{t_{k+1}}]}\phi(s)$ for $t_k\geq T_0$.
Note that it takes at least $h_1(t_k)/h_0$ time for the value of $\phi(t)$ to vary from $h_1(t_k)$ to $0$ or from $0$ to $h_1(t_k)$. By simple calculations, if $h_1(t_k)/h_0\leq d_w/2$, \begin{eqnarray*}\begin{array}{lll}\int_{t_k}^{t_{k+1}}\phi(s)\mathrm{d}s\geq h_1(t_k)^2/h_0/2\end{array}\end{eqnarray*} and if  $h_1(t_k)/h_0\geq d_w/2$, \begin{eqnarray*}\begin{array}{lll}\int_{t_k}^{t_{k+1}}\phi(s)\mathrm{d}s\geq (2h_1(t_k)-h_0d_w/2)d_w/4\geq d_w^2h_0/8.\end{array}\end{eqnarray*} Let $T_0$ be sufficiently large for some $\epsilon<d^2_wh_0/8/(M+1)$ such that $\frac{2\rho}{\sqrt{t}}<\epsilon$
 and {$\|\nabla f_i(x_i(t))\|/\sqrt{q_i(t)}<\epsilon$}
 for all $t\geq T_0$ and all $i$. When $h_1(t_k)>\sqrt{2h_0(M+1)\epsilon}$, \begin{eqnarray*}\begin{array}{lll}\int_{t_k}^{t_{k+1}}\phi(s)\mathrm{d}s\geq \min\{h_1(t_k)^2/h_0/2,d^2_wh_0/8\}\geq
 \epsilon (M+1).\end{array}\end{eqnarray*}  {From Assumption \ref{ass16}, $t_{k+1}-t_k\leq M$. It follows that \begin{eqnarray*}\begin{array}{lll}V_1(t_{k+1})-V_1(t_k)\leq \int_{t_k}^{t_{k+1}}[-\phi(s)+\epsilon]\mathrm{d} s\leq-\epsilon.\end{array}\end{eqnarray*}
 Note that \begin{eqnarray*}\begin{array}{lll}x^T(t)L(t)x(t)\\=\sum_{i=1}^n\sum_{j\in N_i(t)}a_{ij}(t)\|x_i(t)-x_j(t)\|^2\\\geq \sum_{i=1}^n\sum_{j\in N_i(t)}\eta\|x_i(t)-x_j(t)\|^2.\end{array}\end{eqnarray*} When $h_1(t_k)\leq\sqrt{2h_0(M+1)\epsilon}$, $\|x_i(t)-x_j(t)\|\leq \sqrt[4]{2h_0(M+1)\epsilon}/\sqrt{\eta}$ for all $i$ and all $j\in N_i(t)$ and $\|x_i(t)-P_{H_i}(x_i(t))\|\leq \sqrt[4]{2h_0(M+1)\epsilon}$ for all $i$. Thus, $\|\dot{x}_i(t)\|\leq \sqrt[4]{2h_0(M+1)\epsilon}(1+\bar{\eta}n/\sqrt{\eta})+\epsilon$ and $\|x_i(t)-x_i(\bar{t})\|\leq [\sqrt[4]{2h_0(M+1)\epsilon}(1+\bar{\eta}n/\sqrt{\eta})+\epsilon](t_{k+1}-t_k)\leq [\sqrt[4]{2h_0(M+1)\epsilon}(1+\bar{\eta}n/\sqrt{\eta})+\epsilon]M$ for all $t,\bar{t}\in [t_k,t_{k+1}]$.}
  Under Assumption \ref{ass16}, the union of all the graphs during each interval $[{t_k},{t_{k+1}})$, denoted by $\hat{\mathcal{G}}(t_k)$, is strongly connected. In $\hat{\mathcal{G}}(t_k)$, there must exist a directed path between every two agents whose length is no larger than $n$. Thus,
\begin{eqnarray}\label{e23t}\begin{array}{lll}\max_i\{\|x^*(t)-x_i(t)\|\mid t\in [{t_k},{t_{k+1}}]\}\\\leq\max_{i,j}\{\|x_i(t)-x_j(t)\|\mid t\in [{t_k},{t_{k+1}}]\}
\\\leq 2n[\sqrt[4]{2h_0(M+1)\epsilon}(1+\bar{\eta}n/\sqrt{\eta})+\epsilon]M\\+n\sqrt[4]{2h_0(M+1)\epsilon}/\sqrt{\eta}.\end{array}\end{eqnarray}
  From Lemma \ref{lemma7}, when $\lim_{t\rightarrow+\infty}\|x^*(t)-P_{H_i}(x^*(t))\|=0$ for all $i$, $\lim_{t\rightarrow+\infty}\|x^*(t)-P_{H}(x^*(t))\|=0$.
  Let  \begin{eqnarray*}\begin{array}{lll}\phi_1(t)\triangleq\sum_{i=1}^n\|x_i(t)-x^*(t)\|^2\\+\sum_{i=1}^n\|x^*(t)-P_{H_i}(x^*(t))\|^2.\end{array}\end{eqnarray*}
Then there exists a constant $\delta_1>0$ for any ${\epsilon}_1>0$ such that
  $V_1(t)<\epsilon_1$ when $\phi_1(t)<\delta_1$. From Lemma \ref{le1u}, we have
  \begin{eqnarray}\label{e23t1}\begin{array}{lll}\|x^*(t)-P_{H_i}(x^*(t))\|\\=\|x^*(t)-x_i(t)+x_i(t)-P_{H_i}(x_i(t))\\+P_{H_i}(x_i(t))-P_{H_i}(x^*(t))\|\\\leq 2\|x^*(t)-x_i(t)\|+\|x_i(t)-P_{H_i}(x_i(t))\|.\end{array}\end{eqnarray}
Let $\epsilon$ be sufficiently small such that $\max_i\{2\|x^*(t)-x_i(t)\|+\|x_i(t)-P_{H_i}(x_i(t))\|\mid t\in [{t_k},{t_{k+1}}]\}<\sqrt{{\delta}_1/2n}$ when $h_1(t_k)\leq\sqrt{2h_0(M+1)\epsilon}$.
Clearly, $\|x^*(t)-P_{H_i}(x^*(t))\|\leq \sqrt{{\delta}_1/2n}$ and $\|x^*(t)-x_i(t)\|\leq \sqrt{{\delta}_1/2n}$ for all $t\in [{t_k},{t_{k+1}})$ and hence $\phi_1(t)\leq\delta_1$  for all $t\in [{t_k},{t_{k+1}})$ when $h_1(t_k)\leq\sqrt{2h_0(M+1)\epsilon}$. Thus, $V_1(t)<\epsilon_1$  for all $t\in [{t_k},{t_{k+1}})$ when $h_1(t_k)\leq\sqrt{2h_0(M+1)\epsilon}$. Clearly, $h_1(t_k)>\sqrt{2h_0(M+1)\epsilon}$ and hence $V_1(t_{k+1})-V_1(t_k)\leq -\epsilon$ when $V_1(t)\geq\epsilon_1$ for some $t\in [{t_k},{t_{k+1}})$. Thus, there exists a constant $T_1>T_0$ such that $V(t_k)<\epsilon_1$ for all $t_k>T_1$.
 In view of the arbitrariness of $\epsilon_1$, let $\epsilon_1\rightarrow0$. It follows that $\lim_{k\rightarrow+\infty}V_1(t_k)=0$.

From (\ref{e121}),  $\dot{V}_1(t)\leq \frac{2\rho}{\sqrt{t}}$ for all $t\geq T_0$. From Assumption \ref{ass16}, $t_{k+1}-t_k\leq M$. Thus, \begin{eqnarray*}\begin{array}{lll}V_1(t)=V_1(t_k)+\int_{t_k}^t\dot{V}_1(s)\mathrm{d}s\leq V_1(t_k)+\frac{2M\rho}{\sqrt{t_k}}\end{array}\end{eqnarray*} for $t\in [{t_k},{t_{k+1}})$ and hence $\lim_{t\rightarrow+\infty}V_1(t)=0$.
 Thus, $\lim_{t\rightarrow+\infty}\|x^*(t)-P_{H}(x^*(t))\|=\lim_{t\rightarrow+\infty}[x_i(t)-x^*(t)]=0$ for all $i$.
\endproof

\noindent{\textbf{Proof of Theorem \ref{theorem1}:}} Note that $\pi/4\leq\mathrm{arctan}(e^{\|x_i(t)\|})\leq \pi/2$ for all $t$ and all $i$.
There exists a constant $T_0>0$ such that $2\sqrt{t}>\sqrt{q_i(t)}>\frac{\sqrt{t}}{2}$ for all $i$ and all $t\geq T_0$.
From Theorem \ref{lemma1}, all $x_i(t)$ remain in a bounded region, and $\lim_{t\rightarrow+\infty}[x_i(t)-x^*(t)]=\lim_{t\rightarrow+\infty}\|x^*(t)-P_{H}(x^*(t))\|=0$ for all $i$. It is clear that $x^*(t)$ remains in a bounded region. Moreover, $\|x^*(t)-s\|$ is bounded for some $s\in H$ and hence $\|x^*(t)-P_{H}(x^*(t))\|$ is bounded for all $t$. From Lemma \ref{lemma109a1}, $\mathcal{X}$ is a nonempty closed bounded convex set. It follows that $f_i(x^*(t))$ and $f_i(P_{\mathcal{X}}(x^*(t)))$ are both bounded for all $i$ and all $t$. Therefore, together with Lemma \ref{lemma9}, it follows that
 there exists a constant $T_1>T_0$ for any $\epsilon_0>0$ such that \begin{eqnarray*}\begin{array}{lll}\|x^*(t)-x_i(t)\|<\epsilon_0,~~\|x^*(t)-P_{H}(x^*(t))\|<\epsilon_0, \end{array}\end{eqnarray*} \begin{eqnarray*}\begin{array}{lll}{\big|}1-\sqrt{\frac{q_i(t)}{q_j(t)}}{\big|}{\big|}f_i(P_H(x^*(t))){\big|}<\epsilon_0/2\end{array}\end{eqnarray*} and \begin{eqnarray*}\begin{array}{lll}{\big|}1-\sqrt{\frac{q_i(t)}{q_j(t)}}{\big|}{\big|}f_i(P_{\mathcal{X}}(x^*(t))){\big|}<\epsilon_0/2\end{array}\end{eqnarray*} for all $i,j$ and all $t\geq T_1$.

Note that all $f_i(s)$ are differentiable functions and $\lim_{t\rightarrow+\infty}[x_i(t)-x^*(t)]=\lim_{t\rightarrow+\infty}\|x^*(t)-P_{H}(x^*(t))\|=0$. Let $T_1$ be sufficently large such that $|f_i(x^*(t))-f_i(P_{H}(x^*(t)))|<\epsilon_0$ and ${|}f_i(x_i(t))-f_i(P_H(x^*(t))){|}<{\epsilon}_0$  for all $t\geq T_1$.
Consider the Lyapunov function candidate
\begin{eqnarray}\label{e091}{V}_2(t)=\sum_{i=1}^n\|x_i(t)-x^*(t)\|^2+n\|x^*(t)-P_{\chi}(x^*(t))\|^2\end{eqnarray} for all $t\geq T_1$.
Similar to the derivation in (\ref{e121}), we have
\begin{eqnarray*}\begin{array}{lll}\dot{{V}}_2(t)\leq-\phi(t)-2\sum_{i=1}^n(x_i(t)-P_{\chi}(x^*(t)))^T\frac{\nabla
f_i(x_i(t))}{\sqrt{q_i(t)}}.
\end{array}\end{eqnarray*}
Using the convexity of the functions $f_i(s)$, it follows that
\begin{eqnarray*}\begin{array}{lll}\dot{{V}}_2(t)\\
\leq-2\sum_{i=1}^n\frac{f_i(P_H(x^*(t)))-f_i(P_{\chi}(x^*(t)))+f_i(x_i(t))-f_i(P_H(x^*(t)))}{\sqrt{q_i(t)}}\\
\leq-2\sum_{i=1}^n\frac{f_i(P_H(x^*(t)))-f_i(P_{\chi}(x^*(t)))-\epsilon_0}{\sqrt{q_i(t)}}\\
\leq\sum_{i=1}^n\frac{4\epsilon_0}{\sqrt{t}}-2\sum_{i=1}^n\frac{f_i(P_H(x^*(t)))-f_i(P_{\chi}(x^*(t)))}{\sqrt{q_1(t)}}\\
+2\sum_{i=1}^n\frac{f_i(P_H(x^*(t)))-f_i(P_{\chi}(x^*(t)))}{\sqrt{q_1(t)}}{\Big(}1-\sqrt{\frac{q_1(t)}{q_i(t)}}{\Big)}\\
\leq\sum_{i=1}^n\frac{8\epsilon_0}{\sqrt{t}}-\sum_{i=1}^n\frac{f_i(P_H(x^*(t)))-f_i(P_{\chi}(x^*(t)))}{\sqrt{t}}\\
\leq\frac{8n\epsilon_0}{\sqrt{t}}-\sum_{i=1}^n\frac{f_i(x^*(t))-f_i(P_{\chi}(x^*(t)))+f_i(P_H(x^*(t)))-f_i(x^*(t))}{\sqrt{t}}\\
\leq\frac{9n\epsilon_0}{\sqrt{t}}-\sum_{i=1}^n\frac{f_i(x^*(t))-f_i(P_{\chi}(x^*(t)))}{\sqrt{t}}\\
\end{array}\end{eqnarray*}
where the last two inequalities has used the fact that $\frac{\sqrt{t}}{2}<q_i(t)<2\sqrt{t}$ for all $i$ and $t\geq T_0$.

In the following, we use the properties of the point $P_{H}(x^*(t))$ to show the optimal convergence of the system. Define $l_0=\max_{s\in H}\|s-P_\mathcal{X}(s)\|$. Let
\begin{eqnarray*}\begin{array}{lll}E=\{s\in \mathbb{R}^r \mid \|s-P_\mathcal{X}(s)\|\leq l_1\}\end{array}\end{eqnarray*}
 for some constant $0<l_1\leq l_0$. Suppose that $\beta=\min_{s\in H\cap\bar{\partial}E}\sum_{i=1}^n[f_i(s)-f_i(P_\mathcal{X}(s))]$, where $\bar{\partial}E=\{s\in \mathbb{R}^r \mid \|s-P_\mathcal{X}(s)\|= l_1\}$.
 Since $P_\mathcal{X}(s)\in \mathcal{X}$, from the definition of $\mathcal{X}$, we have $\beta>0$. From Lemma \ref{lemma51}, \begin{eqnarray*}\begin{array}{lll}\sum_{i=1}^n[f_i(P_{H}(x^*(t)))-f_i(P_\mathcal{X}(P_{H}(x^*(t))))]\geq\beta\end{array}\end{eqnarray*} for any $P_{H}(x^*(t))\notin E$.
 If $\|x^*(t)-P_{\mathcal{X}}(x^*(t))\|>l_1+\epsilon_0$ for $t\geq T_1$, recalling that $\|x^*(t)-P_{H}(x^*(t))\|<\epsilon_0$ and $|f_i(x^*(t))-f_i(P_{H}(x^*(t)))|<\epsilon_0$ for $t\geq T_1$ and all $i$, it follows that $\|P_{H}(x^*(t))-P_{\mathcal{X}}(x^*(t))\|>l_1$ and \begin{eqnarray*}\begin{array}{lll}\sum_{i=1}^n[f_i(x^*(t))-f_i(P_{\mathcal{X}}(x^*(t)))]>\beta-n\epsilon_0\end{array}\end{eqnarray*} for $t\geq T_1$.
 Let $\epsilon_0$ be sufficiently small such that $2\epsilon_0+10n\epsilon_0<\beta$.
It follows that for any $t\geq T_1$ and any $\|x^*(t)-P_{\mathcal{X}}(x^*(t))\|>l_1+\epsilon_0$, \begin{eqnarray}\label{e1801p}\begin{array}{lll}\dot{{V}}_2(t)\leq -\frac{\beta-n\epsilon_0}{\sqrt{t}}+\frac{9n\epsilon_0}{\sqrt{t}}<-\frac{2\epsilon_0}{\sqrt{t}}.\end{array}\end{eqnarray} Integrating both sides of this inequality from $T_1$ to $t$, we have
 $V_2(t)-{V}_2(T_1)\leq -\epsilon_0(\sqrt{t}-\sqrt{T_1})$.
 This implies that there exists a constant $T_2>T_1$ such that $\|x^*(T_2)-P_{\mathcal{X}}(x^*(T_2))\|=l_1+\epsilon_0$. Since $\|x^*(t)-x_i(t)\|<\epsilon_0$ for all $t\geq T_1$, it follows from (\ref{e091}) that  \begin{eqnarray*}\begin{array}{lll}V_2(T_2)\leq n\epsilon_0^2+n(l_1+\epsilon_0)^2.\end{array}\end{eqnarray*}
 Note that (\ref{e1801p}) holds for any $t\geq T_1$ and any $\|x^*(t)-P_{\mathcal{X}}(x^*(t))\|>l_1+\epsilon_0$.
 Thus, \begin{eqnarray*}\begin{array}{lll}V_2(t)\leq n\epsilon_0^2+n(l_1+\epsilon_0)^2\end{array}\end{eqnarray*} for all $t\geq T_2$.
  In view of the arbitrariness of $l_1$ and $\epsilon_0$, letting $l_1\rightarrow0$ and $\epsilon_0\rightarrow0$,
 we have $\lim_{t\rightarrow+\infty}\|x^*(t)-P_\mathcal{X}(x^*(t))\|=0$. Since $\lim_{t\rightarrow+\infty}[x_i(t)-x^*(t)]=0$, we have $\lim_{t\rightarrow+\infty}\|x_i(t)-P_\mathcal{X}(x^*(t))\|=0$ for all $i$.
That is, all agents reach a consensus and minimize the
team objective function (\ref{gel1}) as $t\rightarrow+\infty$.
\endproof

\noindent{\textbf{Proof of Theorem \ref{theorem14e}:}}
Note that $\pi/4\leq\mathrm{arctan}(e^{\|x_i(k)\|})\leq \pi/2$ for all $k$ and all $i$. There exists a constant $T_0>0$ such that $\frac{kT}{4}<q_i(k)<4kT$ for all $k\geq T_0$.
 Under Assumption \ref{ass2}, $X$ and all $X_i$ are nonempty and bounded. Note from Lemma \ref{lemma109a1} that $\mathcal{X}$ is nonempty and bounded under Assumptions \ref{ass2} and \ref{ass17}. There must exist a closed bounded convex set $Y=\{y\mid \|y-z\|\leq P_2\}\subset \mathbb{R}^r$ for some $z\in \mathcal{X}\subset H$ and some constant $P_2>0$ such that $x_i(T_0)\in Y$,
$X\subset Y$, $X_i\subset Y$ and $\mathcal{X}\subset Y$ for all $i$. Let $Y$ be sufficiently large such that $f_i(v_i(k))-f_i(z)\geq 4\sum_{j=1}^n[f_j(z)-f_j(z_j)]T+4nT$ for all $i$, all $z\in \mathcal{X}$, all $z_j\in X_j$ and all $v_i(k)\notin Y$. 

Consider the Lyapunov function candidate
$V_0(k)=\max_i\|x_i(k)-z\|^2$. Calculating $\|x_i(k+1)-z\|^2$, it follows that
\begin{eqnarray}\label{ep3a}\begin{array}{lll}&\|x_i(k+1)-z\|^2\\
=&\|w_i(k)-z\|^2+\|(P_{H_i}[w_i(k)]-w_i(k))\gamma_i\|^2\\+&2\gamma_i(w_i(k)-z)^T(P_{H_i}[w_i(k)]-w_i(k))\\
\leq& \|w_i(k)-z\|^2+\gamma_i^2\|P_{H_i}[w_i(k)]-w_i(k)\|^2\\+&\gamma_i\|z-P_{H_i}(z)\|^2-\gamma_i\|w_i(k)-P_{H_i}[w_i(k)]\|^2\\
\leq&\|v_i(k)-z-gr_i(k)T\|^2
\end{array}\end{eqnarray}
for all $i$,
where the first inequality has used the convexity of the function $\|s-P_{H_i}(s)\|^2$, and the second inequality has used the fact that $z\in H\subset H_i$ and
$\|z-P_{H_i}(z)\|^2=0$ and the fact that $0<\gamma_i<1$. {When $gr_i(k)=0$, it follows that
\begin{eqnarray}\label{ep2a}\begin{array}{lll}\|v_i(k)-z-gr_i(k)T\|^2
\\\leq \sum\limits_{j\in N_i(k)\cup \{i\}}a_{ij}(k)\|x_j(k)-z\|^2.
\end{array}\end{eqnarray}
When $gr_i(k)\neq0$, $\sqrt{q_i(k)}>\|\nabla
f_i(v_i(k))\|^2$. It follows that
\begin{eqnarray}\label{ep1a}\begin{array}{lll}&\|v_i(k)-z-gr_i(k)T\|^2\\
=&\|v_i(k)-z\|^2+\|\nabla
f_i(v_i(k))T/\sqrt{q_i(k)}\|^2\\-&2(v_i(k)-z)^T\nabla
f_i(v_i(k))T/\sqrt{q_i(k)}\\
\leq &\|v_i(k)-z\|^2+T^2/\sqrt{q_i(k)}\\-&2[f_i(v_i(k))-f_i(z_i)+f_i(z_i)-f_i(z)]T/\sqrt{q_i(k)}\\
\leq &\|v_i(k)-z\|^2+2T^2/\sqrt{kT}\\+&4[f_i(z)-f_i(z_i)]T/\sqrt{kT}\\-&[f_i(v_i(k))-f_i(z_i)]T/\sqrt{kT}
\end{array}\end{eqnarray}
for $k\geq T_0$, where the first inequality has used the convexity of the function $f_i(s)$ and the second inequality has used the fact that $\frac{kT}{4}<q_i(k)<4kT$ for all $k\geq T_0$.} From the definition of $X_i$, $f_i(z)-f_i(z_i)\geq0$ and $f_i(v_i(k))-f_i(z_i)\geq0$. If $v_i(k)\in {Y}$, then $\|x_i(k)-z\|\leq P_2/\eta$ and it follows from (\ref{ep3a}) and (\ref{ep1a}) that
{\begin{eqnarray*}\begin{array}{lll}&\|x_i(k+1)-z\|^2\\\leq& {P}_2^2+2T^2/\sqrt{kT}+4[f_i(z)-f_i(z_i)]T/\sqrt{kT}\\\leq& {P}_2^2+2T^2/\sqrt{T}+4[f_i(z)-f_i(z_i)]\sqrt{T}.\end{array}\end{eqnarray*}}
If $gr_i(k)\neq0$ and $\|x_i(k)-z\|> P_2/\eta$, then $v_i(k)\notin {Y}$ and hence \begin{eqnarray*}\begin{array}{lll}f_{i}(v_{i}(k))-f_{i}(z)\geq 4\sum_{j=1}^n[f_j(z)-f_j(z_j)]+4nT.\end{array}\end{eqnarray*}
It follows from (\ref{ep3a}) and (\ref{ep1a}) that $\|x_i(k+1)-z\|^2\leq V_0(k)$. If $gr_i(k)=0$ for some $i$, it follows from (\ref{ep3a}) and (\ref{ep2a}) that $\|x_i(k+1)-z\|^2\leq V_0(k)$. Summarizing the above analysis, all agents remain in a bounded region. Consequently,
all $\|x_i(k)\|$ and all $\|v_i(k)\|$  are bounded. Since all $f_i(s)$ are differentiable, $\|\nabla f_i(v_i(k))\|$ is bounded for all $i$.
Then there exists a constant $T_1>T_0$ such that $\sqrt{q_i(k)}>\|\nabla
f_i(v_i(k))\|^2$ for all $k\geq T_1$. Thus, \begin{eqnarray*}\begin{array}{lll}gr_i(k)=\frac{\nabla
f_i(v_i(k))}{\sqrt{q_i(k)}}\end{array}\end{eqnarray*} for all $k\geq T_1$.

Now, we prove that $\lim_{k\rightarrow+\infty}\|x^*(k)-P_{H}(x^*(k))\|=\lim_{k\rightarrow+\infty}[x_i(k)-x^*(k)]=0$ for all $i$. Consider the Lyapunov function candidate \begin{eqnarray*}\begin{array}{lll}V_1(k)=V_{11}(k)+V_{12}(k)\end{array}\end{eqnarray*} for $k\geq T_1$, where \begin{eqnarray*}\begin{array}{lll}V_{11}(k)=\sum_{i=1}^n\|x_i(k)-x^*(k)\|^2\end{array}\end{eqnarray*} and \begin{eqnarray*}\begin{array}{lll}V_{12}(k)=n\|x^*(k)-P_{H}(x^*(k))\|^2\end{array}\end{eqnarray*} for $k\geq T_1$. Under Assumption \ref{ass982}, \begin{eqnarray*}\begin{array}{lll}x^*(k+1)\\=\frac{1}{n}\sum_{i=1}^n[v_i(k)-\frac{\nabla
f_i(v_i(k))T}{\sqrt{q_i(k)}}+\gamma_i(P_{H_i}(w_i(k))-w_i(k))]\\=x^*(k)-\frac{1}{n}\sum_{i=1}^n[\frac{\nabla
f_i(v_i(k))T}{\sqrt{q_i(k)}}-\gamma_i(P_{H_i}(w_i(k))-w_i(k))]\end{array}\end{eqnarray*} for $k\geq T_1$.
It follows that
{\begin{eqnarray}\label{eq3s}\begin{array}{lll}&V_{11}(k+1)\\
=&\sum_{i=1}^n[\|w_i(k)+\gamma_i[P_{H_i}(w_i(k))-w_i(k)]\\-&x^*(k+1)\|^2\\
\leq&\sum_{i=1}^n\|w_i(k)-x^*(k)\|^2\\-&\frac{2}{n}\|\sum_{i=1}^n\gamma_i(P_{H_i}(w_i(k))
-w_i(k))\|^2\\+&\frac{1}{n}\|\sum_{i=1}^n\gamma_i(P_{H_i}(w_i(k))-w_i(k))\|^2\\
+&\sum_{i=1}^n\gamma_i^2\|P_{H_i}(w_i(k))-w_i(k)\|^2+\frac{\tilde{c}_1}{\sqrt{k}}+\frac{\tilde{c}_1}{{k}}\\
+&2\sum_{i=1}^n(w_i(k)-x^*(k))^T\gamma_i(P_{H_i}(w_i(k))-w_i(k))\\
-&2\sum_{i=1}^n(w_i(k)-x^*(k))^T\\\times&\frac{1}{n}\sum_{i=1}^n\gamma_i(P_{H_i}(w_i(k))-w_i(k))\\
\leq&\sum_{i=1}^n\|w_i(k)-x^*(k)\|^2\\-&\frac{1}{n}\|\sum_{i=1}^n\gamma_i(P_{H_i}(w_i(k))
-w_i(k))\|^2\\+&\sum_{i=1}^n\gamma_i^2\|P_{H_i}(w_i(k))-w_i(k)\|^2\\
+&2\sum_{i=1}^n(w_i(k)-x^*(k))^T\gamma_i(P_{H_i}(w_i(k))-w_i(k))\\
+&\frac{c_1}{\sqrt{k}}+\frac{c_1}{{k}}\\
\end{array}\end{eqnarray}}
\hspace{-0.2cm}for all $k\geq T_1$ and two constants $0<\tilde{c}_1<c_1$, where \begin{eqnarray*}\hspace{-0.2cm}\begin{array}{lll}&|\sum_{i=1}^n(w_i(k)-x^*(k))^T\frac{1}{n}\sum_{i=1}^n\gamma_i(P_{H_i}(w_i(k))-w_i(k))|\\
&=|\sum_{i=1}^n(x^*(k)-x^*(k)-\frac{1}{n}\sum_{i=1}^n\frac{\nabla
f_i(v_i(k))T}{\sqrt{q_i(k)}})^T\\&\times\frac{1}{n}\sum_{i=1}^n\gamma_i(P_{H_i}(w_i(k))-w_i(k))|\leq \frac{c_1-\tilde{c}_1}{2\sqrt{k}}\end{array}\end{eqnarray*} under Assumption \ref{ass982} and the terms such as $\frac{c_1}{\sqrt{k}}$ and $\frac{c_1}{{k}}$ can be obtained based on the fact that the variables such as $\|P_{H_i}(w_i(k))\|$, $\|w_i(k)\|$ and $\|\nabla
f_i(v_i(k))\|$ are all bounded and the fact that $\frac{kT}{4}<q_i(k)<4kT$ for all $k\geq T_0$. Also,
{\begin{eqnarray}\label{eq3s1}\hspace{-0.1cm}\begin{array}{lll}&V_{12}(k+1)
\leq n\|x^*(k+1)-P_{H}(x^*(k))\|^2\\
\leq&n\|x^*(k)-P_{H}(x^*(k))-\frac{1}{n}\sum_{i=1}^n[\frac{\nabla
f_i(v_i(k))}{\sqrt{q_i(k)}}\\-&\gamma_i(P_{H_i}(w_i(k))-w_i(k))]\|^2\\
\leq&n\|x^*(k)-P_{H}(x^*(k))\|^2+2\sum_{i=1}^n[x^*(k)-w_i(k)\\+&w_i(k)-P_{H}(x^*(k))]^T\gamma_i[P_{H_i}(w_i(k))-w_i(k))]\\
+&\frac{c_2}{\sqrt{k}}+\frac{c_2}{{k}}+\frac{1}{n}\|\sum_{i=1}^n\gamma_i(P_{H_i}(w_i(k))-w_i(k))\|^2\\
\leq&n\|x^*(k)-P_{H}(x^*(k))\|^2\\+&2\sum_{i=1}^n[x^*(k)-w_i(k)]^T\\
\times&\gamma_i[P_{H_i}(w_i(k))-w_i(k))]+\frac{c_2}{\sqrt{k}}+\frac{c_2}{{k}}\\
-&\sum_{i=1}^n\gamma_i\|P_{H_i}(w_i(k))-w_i(k))\|^2\\
+&\frac{1}{n}\|\sum_{i=1}^n\gamma_i(P_{H_i}(w_i(k))-w_i(k))\|^2
\end{array}\end{eqnarray}}
\noindent\hspace{-0.2cm}for all $k\geq T_1$ and some constant $c_2>0$,
where the last inequality has used the convexity of $\frac{1}{2}\|P_{H_i}(w_i(k))-w_i(k)\|^2$, i.e., \begin{eqnarray*}\begin{array}{lll}[w_i(k)-P_{H}(x^*(k))]^T[P_{H_i}(w_i(k))-w_i(k))]\\\leq \frac{1}{2}\|P_{H_i}(P_{H}(x^*(k)))-P_{H}(x^*(k))\|^2\\-\frac{1}{2}\|P_{H_i}(w_i(k))-w_i(k)\|^2\leq -\frac{1}{2}\|P_{H_i}(w_i(k))-w_i(k)\|^2.\end{array}\end{eqnarray*}

 Under Assumption \ref{ass982}, we have $\mathcal{A}(k)\mathbf{1}=\mathbf{1}$ and hence $\mathcal{A}(k)\mathbf{1}x^*(k)=\mathbf{1}x^*(k)$. It follows that \begin{eqnarray*}\begin{array}{lll}\sum_{i=1}^n\|v_i(k)-x^*(k)\|^2\\=\|\mathcal{A}(k)x(k)-\mathbf{1}x^*(k)\|^2=\|\mathcal{A}(k)[x(k)-\mathbf{1}x^*(k)]\|^2.\end{array}\end{eqnarray*}
 Note from the form of $w_i(k)$ that $\sum_{i=1}^n\|w_i(k)-x^*(k)\|^2\leq \sum_{i=1}^n\|v_i(k)-x^*(k)\|^2+\frac{c_3}{\sqrt{k}}+\frac{c_3}{{k}}$
 for some constant $c_3>0$ and $k\geq T_1$. Together with (\ref{eq3s}) and (\ref{eq3s1}), we have that \begin{eqnarray*}\begin{array}{lll}V_1(k+1)-V_1(k)\leq -\bar{\phi}(k) +\frac{c_4}{\sqrt{k}}+\frac{c_4}{{k}},\end{array}\end{eqnarray*} for some constant $c_4>0$ and $k\geq T_1$ where \begin{eqnarray*}\begin{array}{lll}\bar{\phi}(k)=[x(k)-\mathbf{1}x^*(k)]^T[I-\mathcal{A}(k)^T\mathcal{A}(k)][x(k)-\mathbf{1}x^*(k)]\\+ \sum_{i=1}^n(1-\gamma_i)\gamma_i\|P_{H_i}(w_i(k))-w_i(k))\|^2.\end{array}\end{eqnarray*} By repeated calculations,
 \begin{eqnarray*}\begin{array}{lll}V_1(k_{m+1})-V_1(k_m)\leq  -\sum_{k=k_m}^{{k_{m+1}-1}}\bar{\phi}(k)+\frac{c_5}{\sqrt{k_m}}+\frac{c_5}{{k_m}}\end{array}\end{eqnarray*} for some constant $c_5>0$.
 Let $T_1$ be sufficiently large for any $\epsilon>0$ such that $\frac{c_5}{\sqrt{k_m}}+\frac{c_5}{{k_m}}<\epsilon$ for all $k\geq T_1$ and $h_1(k_m)=\max_{k\in [{k_m},{k_{m+1}})}\bar{\phi}(k)$ for $k_m\geq T_1$.
  When $h_1(k_m)>2\epsilon$ for $k_m\geq T_1$,
 \begin{eqnarray*}\begin{array}{lll}V_1(k_{m+1})-V_1(k_m)<-\epsilon.\end{array}\end{eqnarray*}

 Now, we need to consider the upper bound of $\|x_i(k)-x^*(k)\|$ when $h_1(k_m)\leq 2\epsilon$ for $k_m\geq T_1$. Under Assumption \ref{ass982}, the graph $\mathcal{G}(k)$ is balanced for all $k$. Rearranging the agent indices, it can be obtained that $\mathcal{A}(k)$ can be denoted as $\mathcal{A}(k)=\mathrm{diag}\{\mathcal{A}^1(k),\cdots,\mathcal{A}^h(k)\}$, where each $\mathcal{A}^i(k)$ corresponds to a strongly connected component of the agents. As a result, \begin{eqnarray*}\begin{array}{lll}\mathcal{A}(k)^T\mathcal{A}(k)=\mathrm{diag}\{\mathcal{A}^1(k)^T\mathcal{A}^1(k),\cdots,\mathcal{A}^h(k)^T\mathcal{A}^h(k)\}.\end{array}\end{eqnarray*} Note that {$\sum_{i=1}^na_{ij}(k)=\sum_{j=1}^n a_{ji}(k)=1$} under Assumption \ref{ass982}. It follows that each row sum of $\mathcal{A}^i(k)^T\mathcal{A}^i(k)$ is $1$ and hence $[I-\mathcal{A}(k)^T\mathcal{A}(k)]\mathbf{1}=0$.
 Note that $\mathcal{A}^i(k)^T\mathcal{A}^i(k)$ is symmetric and each of its off-diagonal entries is nonnegative. $I-\mathcal{A}(k)^T\mathcal{A}(k)$ can be regarded as the Laplacian of a certain undirected graph, denoted by $\bar{\mathcal{G}}(k)$, and each $I-\mathcal{A}^i(k)^T\mathcal{A}^i(k)$ corresponds to each connected component of $\bar{\mathcal{G}}(k)$, which in turn corresponds to the strongly connected component of ${\mathcal{G}}(k)$.  Under Assumption \ref{ass16o}, the union of all $\mathcal{G}(k)$ during each interval $[k_j,k_{j+1})$
is strongly connected. Hence, the union of all $\bar{\mathcal{G}}(k)$ during each interval $[k_j,k_{j+1})$ is connected.
 Note that each nonzero entry of $\mathcal{A}(k)$ is no smaller than $\eta$ under Assumption \ref{ass982} and hence each nonzero entry of $\mathcal{A}(k)^T\mathcal{A}(k)$ is no smaller than $\eta^2$. Thus, \begin{eqnarray*}\begin{array}{lll}x(k)^T([\mathcal{A}(k)^T\mathcal{A}(k)-I_n]x(k)\\\leq -\sum_{i=1}^n\sum_{j\in N_i(k)}\eta^2\|x_i(k)-x_j(k)\|^2.\end{array}\end{eqnarray*} By some calculations similar to (\ref{e23t}) and (\ref{e23t1}), there exists a constant $c_6>0$ such that $\|x_i(k)-x^*(k)\|\leq c_6 (\sqrt{\epsilon}+\epsilon)$ and $\|x^*(k)-P_{H_i}(x^*(k))\|\leq c_6 (\sqrt{\epsilon}+\epsilon)$ for all $i$ and all $k_m\leq k<k_{m+1}-1$ when $h_1(k_m)\leq 2\epsilon$ for $k_m>T_1$. Let \begin{eqnarray*}\begin{array}{lll}\bar{\phi}_1(k)&\triangleq\sum_{i=1}^n\|x_i(k)-x^*(k)\|^2\\&+\sum_{i=1}^n\|x^*(k)-P_{H_i}(x^*(k))\|^2. \end{array}\end{eqnarray*}
 From Lemma \ref{lemma7}, when $\lim_{t\rightarrow+\infty}\|x^*(k)-P_{H_i}(x^*(k))\|=0$ for all $i$, $\lim_{t\rightarrow+\infty}\|x^*(k)-P_{H}(x^*(k))\|=0$.
  There exists a constant $\delta_1>0$ for any $\epsilon_1>0$ such that $V_1(k)<\epsilon_1$ when $\bar{\phi}_1(k)<\delta_1$. Let $\epsilon$ be sufficiently small such that $\bar{\phi}_1(k)\leq{\delta}_1$ and hence $V_1(k)<\epsilon_1$ for all $k_m\leq k<k_{m+1}$ when $h_1(k_m)\leq 2\epsilon$ for $k_m>T_1$. Note that $h_1(k_m)>2\epsilon$ when $V_1(k)>\epsilon_1$ for some $T_1<k_{m}\leq k<k_{m+1}$ and recall that
when $h_1(k_m)>2\epsilon$,
 $V_1(k_{m+1})-V_1(k_m)<-\epsilon$. Thus, there exists a constant $T_1>T_0$ such that
 $V_1(k_m)\leq \epsilon_1$ for all $k_m\geq T_1$. In view of the arbitrariness of $\epsilon_1$, let $\epsilon_1\rightarrow0$. It follows that $\lim_{k_m\rightarrow+\infty}V_1(k_m)=0$. Since \begin{eqnarray*}\begin{array}{lll}V_1(k+1)-V_1(k)\leq -\bar{\phi}(k)+\frac{c_4}{\sqrt{k}}+\frac{c_4}{{k}}\leq \frac{c_4}{\sqrt{k}}+\frac{c_4}{{k}}\end{array}\end{eqnarray*} for all $k>T_1$ and $k_{m+1}-k_m\leq M$, it follows that $\lim_{k\rightarrow+\infty}V_1(k)=\lim_{k_m\rightarrow+\infty}V_1(k_m)=0$.
 Thus, $\lim_{k\rightarrow+\infty}\|x^*(k)-P_{H}(x^*(k))\|=\lim_{k\rightarrow+\infty}[x_i(k)-x^*(k)]=0$ for all $i$. The rest proof is very similar to that of Theorem \ref{theorem1} and hence omitted.
\endproof
\noindent{\textbf{Proof of Theorem \ref{theorem2}:}} Since $x_i(0)\in H_i$ and $\gamma_i=1$ for all $i$, from (\ref{eq53371}), $x_i(k)\in H_i$ for all $i$ and all $k$.
 By the same approach of the proof of Theorem \ref{theorem14e}, this theorem can be proved. However, it should be noted that $\|x^*(k)-P_{H_i}(x^*(k))\|\leq \|x_i(k)-x^*(k)\|\leq c_6 (\sqrt{\epsilon}+\epsilon)$ when it has been proved that $\|x_i(k)-x^*(k)\|\leq c_6 (\sqrt{\epsilon}+\epsilon)$ for all $i$ and two certain constants $c_6>0$ and $\epsilon>0$.
\endproof

\end{document}